\documentclass[12pt]{amsart}
\usepackage{amssymb}
\usepackage{graphicx}

\newtheorem{theorem}{Theorem}

\newtheorem{lemma}[theorem]{Lemma}

\newtheorem{cor}[theorem]{Corollary}

\newtheorem{exam}[theorem]{Example}

\newtheorem{remark}[theorem]{Remark}

\newtheorem{defi}[theorem]{Definition}

\newcommand{\eqa}{\begin{eqnarray}}
\newcommand{\eeqa}{\end{eqnarray}}
\newcommand{\beq}{\begin{equation}}
\newcommand{\eeq}{\end{equation}}
\newcommand{\nn}{\nonumber}

\newcommand{\pal}{\partial}

\newcommand{\pf}{\noindent{\it Proof \ }}

\newcommand{\ub}{{\bf u}}
\newcommand{\epf}{$\quad$\hfill
\raisebox{0.11truecm}{\fbox{}}\par\vskip0.4truecm}

\begin{document}

\title[Universality in Hamiltonian PDEs]
{On universality of critical behaviour in Hamiltonian PDEs}

\author[B. Dubrovin]{Boris Dubrovin}

\address{SISSA\\ Via Beirut 2--4\\ 34014 Trieste\\ Italy
\\ and 
Steklov Math. Institute\\ Moscow}

\email{dubrovin@sissa.it}

\maketitle

\begin{flushright}
\parbox{9cm}{
\begin{center}
{\it Dedicated
to Sergei Petrovich Novikov \\
on the occasion of his 70th birthday.}
\end{center}
}
\end{flushright}

%%%%%%%%%%%%%%%%%%%%%%%%%%%%%%%%%%%%%%%%%%%%%
\begin{abstract}
Our main goal is the comparative study of singularities of solutions to the systems of first order quasilinear PDEs and their perturbations containing higher derivatives. The study is focused on the subclass of Hamiltonian PDEs with one spatial dimension. For the systems of order one or two we describe the local structure of singularities of a generic solution to the unperturbed system near the point of ``gradient catastrophe" in terms of standard objects of the classical singularity theory; we argue that their perturbed companions must be given by certain special solutions of Painlev\'e equations and their generalizations. 
\end{abstract}
%%%%%%%%%%%%%%%%%%%%%%%%%%%%%%%%%%%%%%%%%%%%%

\tableofcontents

\section{Introduction}
In \cite{DZ} we developed a perturbative approach to classification of systems of integrable Hamiltonian partial differential equations (PDEs). The systems under consideration are written in the form of {\it long wave expansion} (also called {\it small dispersion expansion})
\beq\label{pert1}
\ub_t=A(\ub )\, \ub_x +B_2(\ub; \ub_x, \ub_{xx}) +B_3 (\ub; \ub_x, \ub_{xx}, \ub_{xxx})+\dots .
\eeq
Here
$$
\ub=\left( u^1(x,t), \dots, u^n(x,t)\right)^{\rm T}
$$
is the unknown vector function, the entries of the coefficient matrix $A(\ub)$ are smooth on some domain ${\mathcal D}\subset {\mathbb R}^n$; the characteristic roots $\lambda_1(\ub)$, \dots, $\lambda_n(\ub)$ will be assumed to be pairiwise distinct 
\beq\label{hyp1}
\det \left( A(\ub)-\lambda\cdot 1\right) =0, \quad \lambda_i(\ub) \neq \lambda_j(\ub)\quad \mbox{for}\quad i\neq j, \quad \forall \,\ub\in{\mathcal D}.
\eeq
The terms $B_2$, $B_3$, \dots, $B_k$, \dots of the expansion\footnote{In order to uniformize the notations the leading term could be redenoted as
$$
A(\ub)\, \ub_x =: B_1(\ub; \ub_x).
$$
However we will not follow this system of notations. We emphasize that the terms of degree zero (thus, depending only on $\ub$ but not on the jets) must not appear in the rhs of \eqref{pert1}.}
in the rhs of \eqref{pert1} are polynomials in the jet coordinates
$$
\ub_x=\left(u^1_x, \dots, u^n_x\right), \quad \ub_{xx}=\left( u^1_{xx}, \dots, u^n_{xx}\right), \quad\dots
$$
graded homogeneous of the degrees 2, 3, \dots, $k$, \dots,
\eqa\label{pert2}
&&
\deg B_k\left(\ub; \ub_x, \dots, \ub^{(k)}\right) =k
\nn\\
&&
\\
&&
\deg \frac{\pal ^m u^i}{\pal x^m} =m, \quad m>0, \quad \deg u^i=0, \quad i=1, \dots, n.
\nn
\eeqa
The coefficients of these polynomials are smooth functions on the same domain ${\mathcal D}$. 

The system \eqref{pert1} can be considered as a perturbation of the first order quasilinear system
\beq\label{unpert}
\ub_t=A(\ub)\, \ub_x
\eeq
when considering slow varying solutions, i.e., working with vector functions $\ub(x)$ changing by 1 on the spatial scale of order $L\gg1$. The magnitude $|\ub (x)|$ will {\it not} be assumed to be small; the natural small parameter 
\beq\label{small1}
\epsilon=\frac1{L}
\eeq
appears in the estimates for the derivatives
\beq\label{small2}
\ub_x \sim \epsilon, \quad \ub_{xx}\sim \epsilon^2, \dots, \ub^{(k)}\sim \epsilon^k, \dots.
\eeq
It is convenient to introduce {\it slow variables} by rescaling
\beq\label{slow}
x\mapsto \epsilon\, x, \quad t\mapsto \epsilon\, t.
\eeq
The system \eqref{pert1} rewrites in the form
\beq\label{pert0}
\ub_t=A(\ub )\, \ub_x +\epsilon\,B_2(\ub; \ub_x, \ub_{xx}) +\epsilon^2 B_3 (\ub; \ub_x, \ub_{xx}, \ub_{xxx})+\dots .
\eeq
In these notations the leading term \eqref{unpert} results from  the limit $\epsilon\to 0$.

For sufficiently small $\epsilon$ one expects to see no major differences in the behaviour of solutions to the perturbed and unperturbed equations \eqref{pert0} and \eqref{unpert} within the regions where the $x$-derivatives
are bounded. However the differences become quite serious near the critical point (also called the point of {\it gradient catastrophe}) where derivatives of the solution to the unperturbed equation tend to infinity.

Although the case of small viscosity perturbations has been well studied and understood (see \cite{bressan} and references therein), the critical behaviour of solutions to {\it general conservative perturbations}
\eqref{pert0} to our best knowledge has not been investigated in a systematic way before the author's paper \cite{du2}
(see the papers \cite{gp, hl, ks, ll, llv} for the study of various particular cases).

The main goal of this paper is to extend the technique of \cite{du2} to Hamiltonian perturbations of systems of PDEs.
A serious difference from the scalar case treated in \cite{du2} can be seen already in the case of second order systems considered in the present paper. For the scalar case {\it any} Hamiltonian perturbation preserves integrability up to a high order of the small parameter $\epsilon$. For systems this is not the case: for generic perturbation of a second order quasilinear system\footnote{Recall that any second order quasilinear PDE \eqref{unpert} can be integrated by the method of characteristics.} integrability is destroyed already at the order $\epsilon^2$.

Due to this difference we begin our study of critical behaviour in perturbed systems of Hamiltonian PDEs with considering {\it integrable perturbations} only.  Our idea is that the local structure of the critical behaviour should depend only on the order $\epsilon^4$ perturbation near a weak singularity of a hyperbolic system on the order $\epsilon^2$ perturbation near a weak singularity of an elliptic system. (We postpone to a separate publication the study of singularities on the boundary between domains of ellipticity and hyperbolicity.) Because of this it suffices to have integrability only at low orders of the perturbative expansion. We give a precise meaning to such an ``approximate integrability", explain the basic tools involved
(extension of infinitesimal symmetries, $D$-operator, quasitriviality and string equation), and formulate 
 the Universality Conjectures
about the behaviour of a generic solution to the  perturbed Hamiltonian system near the point of gradient catastrophe of the unperturbed solution. We argue that, up to shifts, Galilean transformations and rescalings this behaviour essentially {\it does not} depend on the choice of solution  (provided certain genericity assumptions hold true). Moreover, this behaviour near the critical point is given by the following formulae:

\noindent $\bullet$ near a critical point within the domain of hyperbolicity the solution, after an affine transformation of the independent variables $(x,t)\to (x_+, x_-)$ and a nonlinear change of the dependent variables $(u,v)\mapsto (r_+, r_-)$ the solution has the following asymptotic representation
\eqa\label{univer}
&&
r_+ \simeq r_+^0 + c\, x_++\alpha_+ \epsilon^{4/7} U''\left( a\, \epsilon^{-6/7} x_-; b\, \epsilon^{-4/7} x_+\right) 
+{\mathcal O}\left( \epsilon^{6/7}\right)
\nn\\
&&
\\
&&
r_- \simeq r_-^0 +\alpha_-\,\epsilon^{2/7} U \left( a\, \epsilon^{-6/7} x_-; b\, \epsilon^{-4/7} x_+\right) +{\mathcal O}\left( \epsilon^{4/7}\right)
\nn
\eeqa
where $U=U(X; T)$ is a particular real {\it smooth for all $X\in \mathbb R$} solution to the fourth order analogue of the Painlev\'e-I equation
\beq\label{ode0}
X=T\, U -\left[ \frac16 U^3 +\frac1{24} ( {U'}^2 + 2 U\, U'' ) +\frac1{240} U^{IV}\right], \quad U'=\frac{dU}{dX} \quad \mbox{ etc.}
\eeq
depending on the real parameter $T$ (the so-called $P_I^2$ equation). Here $a$, $b$, $c$, $\alpha_\pm$ are some constants that depend on the choice of the equation and the solution. For the original variables $u$, $v$ the above asymptotic behaviour implies
\eqa\label{clvl}
&&
u\simeq u_{\rm c} + \lambda\, \epsilon^{2/7} U \left( a\, \epsilon^{-6/7} x_-; b\, \epsilon^{-4/7} x_+\right) +{\mathcal O}\left( \epsilon^{4/7}\right)
\nn\\
&&
\\
&&
v\simeq v_{\rm c} + \mu\, \epsilon^{2/7} U \left( a\, \epsilon^{-6/7} x_-; b\, \epsilon^{-4/7} x_+\right) +{\mathcal O}\left( \epsilon^{4/7}\right)
\nn
\eeqa
with some constants $\lambda$, $\mu$.

The solution $U(X,T)$ was completely characterized by T.Claeys and M.Vanlessen in \cite{cl1} via the isomonodromy deformations technique (see \cite{kapaev}) where a long standing problem of existence of a pole-free real solution to $P_I^2$ has been settled\footnote{The equation (\ref{ode0}) appeared in \cite{bmp} (for the particular value  of the parameter $T=0$) in the study of the double scaling limit for the matrix model with the multicritical index $m=3$. It was observed that the generic solution to (\ref{ode0}) blows up at some point of real line; the conjecture about existence of a unique smooth solution has been formulated. Our original conjecture formulated in \cite{du2} says that the solution to \eqref{ode0} smooth on the entire real axis exists and is {\it unique} for any real $T$. The uniqueness part of this conjecture remains open.}. In the very recent paper \cite{cg} by T.Claeys and T.Grava the Universality Conjecture formulated in \cite{du2} was rigorously established, with the help of Deift - Zhou steepest descent method \cite{dz}, for the case of critical behaviour of solutions to the KdV equation with analytic initial data. For the systems the somewhat weaker asymptotics \eqref{clvl} was proved in \cite{cl2} for a particular class of solution to the hierarchy of Toda equations with analytic initial data. These are the solutions appearing in the study of Hermitean matrix integrals. They depend on one arbitrary analytic function that determines the probability distribution; the general solution to Toda hierarchy equations should depend on two arbitrary functions. The critical behaviour \eqref{clvl} appears, according to \cite{cl2}, in the large $N$ asymptotic behaviour of orthogonal polynomials near the so-called singular edge points.
The small dispersion parameter in this case is $\epsilon=1/N$.

\noindent $\bullet$ near a critical point within the domain of ellipticity a suitable complex combination $w=w(u,v)$ of the solution $u(x,t)$, $v(x,t)$ as function of an affine complex combination $z $ of the independent variables behaves as follows
\eqa\label{e-univer}
&&
w \simeq w^0 +\alpha\, \epsilon^{2/5} W_0\left(  \epsilon^{-4/5}\, z \right) +{\mathcal O}\left( \epsilon^{4/5}\right)
\nn\\
&&
\\
&&
z= a\, x + b\, t + z_0, \quad a, \, b \in \mathbb C.
\nn
\eeqa
Here $W_0=W_0(Z)$ is the so-called {\it tritronqu\'ee} solution to the Painlev\'e-I equation ($P_I$)
$$
W''=6\, W^2 -Z.
$$
This solution discovered by Boutroux in \cite{bo} is uniquely characterized by the following asymptotic behaviour
$$
W(Z) \sim -\sqrt{\frac{Z}6}, \quad |Z|\to\infty, \quad |\arg Z| <\frac{4\pi}5.
$$
In particular it has no poles for sufficiently large $|Z|$ in the above sector.

In \cite{dgk} it was conjectured that the {\it tritronqu\'ee} solution has {\it no poles} inside the sector
\beq\label{sec}
|\arg Z| <\frac{4\pi}5.
\eeq
Although a rigorous proof of this conjecture remains an open problem, numerical evidences found in \cite{dgk} support the conjectur.

Importance of absence of poles in the sector \eqref{sec} has the following explanation. It turns out that the complex parameters $a$ and $b$ in \eqref{e-univer} are such that for sufficiently small $|t|$ and $\epsilon>0$ the argument 
$$
Z=\frac{a\, x+b\, t + z_0}{\epsilon^{4/5}}
$$
of the {\it tritronqu\'ee} solution belongs to the sector \eqref{sec}
for all real $x$. So the Universality Conjecture \eqref{e-univer} makes sense only if no poles occur within the sector \eqref{sec}.

The Universality Conjecture \eqref{e-univer} appeared first in the recent paper \cite{dgk} for the case of focusing nonlinear Schr\"odinger (NLS) equation. In this paper we will explain that the same critical behaviour arises also for other Hamiltonian perturbations of second order quasilinear PDEs of {\it elliptic} type.
 
The paper is organized as follows. In Section 2 we recall some basic facts about second order quasilinear PDEs, their Hamiltonian structures, solutions, conserved quantities and symmetries. In Section 3 we give a very simple local classification of weak singularities of solutions, for both hyperbolic and elliptic cases. An attempt to extend this local classification of singularities to quasilinear hyperbolic systems of order higher than 2 is also given; only the case of integrable hyperbolic PDEs (the so-called {\it semi-hamiltonian} systems in the terminology of \cite {tsarev}), or {\it syst\`emes rich} according to \cite{serre}  is considered. In Section 4 the main tools for the study of Hamiltonian perturbations are developed. Applications of these techniques to the study of concrete examples is given in Section 5.

\vskip 0.5truecm
\noindent{\bf Acknowledgments.} This work is
partially supported by European Science Foundation Programme ``Methods of
Integrable Systems, Geometry, Applied Mathematics" (MISGAM), Marie Curie RTN ``European Network in Geometry, Mathematical Physics and Applications"  (ENIGMA), 
and by Italian Ministry of Universities and Researches (MUR) research grant PRIN 2006
``Geometric methods in the theory of nonlinear waves and their applications". 

\setcounter{equation}{0}
\setcounter{theorem}{0}
\section{Nonlinear wave equation: Hamiltonian structure, first integrals and solutions}\par

We now proceed to considering the two-component systems of Hamiltonian PDEs. In this section we start with the first order quasilinear systems \eqref{unpert}. Moreover, we restrict ourselves with considering an important subclass of such systems associated with nonlinear wave equation
\beq\label{wave1}
u_{tt} -\pal_x^2 P(u)=0
\eeq
for a given smooth function $P(u)$. The equation \eqref{wave1}
is linear for a quadratic function $P(u)$; we assume therefore that
$$
P'''(u)\neq0.
$$

\subsection{Main examples of nonlinear wave equations}\label{sect21}\par

The equation \eqref{wave1} of course is of interest {\sl per se}; it also arises in the study of dispersionless limits of various PDEs of higher order. In particular, for $P(u)=-\frac16\, u^3$
one obtains the dispersionless limit
\beq\label{bous0}
u_{tt} + \left( u\, u_x\right)_x=0
\eeq
of Boussinesq equation
\beq\label{bous}
u_{tt} + \left( u\, u_x\right)_x+u_{xxxx}=0.
\eeq
For $P(u)=e^u$ \eqref{wave1} reduces to the long-wave limit
\beq\label{toda0}
u_{tt} =\pal_x^2 e^u
\eeq
of Toda equations
\beq\label{toda}
\ddot u_n =e^{u_{n+1}-u_n}-e^{u_n-u_{n-1}}
\eeq
More generally, one-dimensional system of particles with neighboring interaction
\beq\label{fpu1}
H=\sum \frac12 p_n^2 +P(q_n-q_{n-1})
\eeq
with the potential $P(u)$ (generalized Fermi - Pasta - Ulam system) yields, after interpolation
\eqa\label{interpol}
&&
q_n(t)-q_{n-1}(t) =w(\epsilon\, n, \epsilon\, t)
\nn\\
&&
\\
&&
p_n(t) = v(\epsilon\, n, \epsilon\, t)
\nn
\eeqa
%the following system
%\eqa
%&&
%u_t = \frac{1-e^{-\epsilon\, \pal_x}}{\epsilon}\, w
%\nn\\
%&&
%w_t=\frac{e^{\epsilon\, \pal_x}-1}{\epsilon}\, P'(u).
%\nn
%\eeqa
and substitution
\beq\label{uw}
u=\frac{\epsilon\, \pal_x}{1-e^{-\epsilon\,\pal_x}}\, w
\eeq
the following system
\eqa\label{fpu2}
&&
u_t=v_x
\nn\\
&&
\\
&&
v_t= \epsilon^{-1}\left[P'\left(\frac{e^{\epsilon\,\pal_x} -1}{\epsilon\,\pal_x}u\right)-P'\left(\frac{1-e^{-\epsilon\,\pal_x} }{\epsilon\,\pal_x}u\right)\right]
\nn\\
&&
\nn\\
&&
\quad= \pal_x P'(u) +\frac{\epsilon^2}{24}\left[ 2\, P''(u)\, u_{xxx} + 4\, P'''(u)\, u_x u_{xx} + P^{IV}(u)\, u_x^3\right] +{\mathcal O}(\epsilon^4).
\nn
\eeqa
%where
%$$
%A_\epsilon= \left(\frac{\sinh \frac{\epsilon\, \pal_x}2}{\frac{\epsilon\, \pal_x}2}\right)^2=1+\frac1{12} \left(\epsilon\, \pal_x\right)^2 +\frac1{360}\left(\epsilon\, \pal_x\right)^4+\dots .
%$$
The above formulae are understood as formal power series in $\epsilon$:
$$
w=\frac{e^{\epsilon\,\pal_x} -1}{\epsilon\,\pal_x}\,u =\frac1{\epsilon}\,\int_x^{x+\epsilon} u(s)\, ds= u +\sum_{k\geq 1} \frac{\epsilon^k}{(k+1)!} \, u^{(k)},
$$
\beq\label{bernoul}
u= w+\frac12\,\epsilon\, w' + \sum_{k > 1} \frac{B_k}{k!} \epsilon^k w^{(k)},
\eeq
 $B_k$ are the Bernoulli numbers.
In the limit $\epsilon\to 0$ \eqref{fpu2} reduces to the nonlinear wave equation written as a system
\eqa\label{wave2}
&&
u_t=v_x
\nn\\
&&
\\
&&
v_t=\pal_x P'(u).
\nn
\eeqa

Further examples will be considered in Section \ref{sect5}. 

\subsection{Hamiltonian formulation of nonlinear wave equation}\par

The system \eqref{wave2} can be written in the Hamiltonian form
\eqa\label{wave3}
&&
u_t=\pal_x\frac{\delta H}{\delta v(x)}
\nn\\
&&
\\
&&
v_t=\pal_x \frac{\delta H}{\delta u(x)}
\nn
\eeqa
with the Hamiltonian
\beq\label{hwave}
H=\int \left[ \frac12 v^2 +P(u)\right]\,dx.
\eeq
The associated Poisson bracket reads
\beq\label{pb}
\{ u(x), v(y)\}=\delta'(x-y)
\eeq
The Poisson bracket of two local functionals
$$
F=\int f(u, v; u_x, v_x, \dots)\, dx, \quad G=\int g(u,v; u_x, v_x, \dots)\, dx
$$
is again a local functional given by
\beq\label{pb1}
\{ F, G\} =\int \left[ \frac{\delta F}{\delta u(x)} \,\pal_x 
\frac{\delta G}{\delta v(x)}+\frac{\delta F}{\delta v(x)}\, \pal_x 
\frac{\delta G}{\delta u(x)}\right]\, dx
\eeq
where
\eqa
&&
\frac{\delta F}{\delta u(x)} ={\rm E}_u f, \quad \frac{\delta G}{\delta u(x)} ={\rm E}_u g
\nn\\
&&
\frac{\delta F}{\delta v(x)} ={\rm E}_v f, \quad \frac{\delta G}{\delta v(x)} ={\rm E}_v g
\nn
\eeqa
are the variational derivatives, 
\eqa
&&
{\rm E}_u:=
\frac{\pal }{\pal u} -\pal_x \frac{\pal }{\pal u_x}+\pal_x^2 \frac{\pal }{\pal u_{xx}}-\dots
\nn\\
&&
\\
&&
{\rm E}_v:=
\frac{\pal }{\pal v} -\pal_x \frac{\pal }{\pal v_x}+\pal_x^2 \frac{\pal }{\pal v_{xx}}-\dots
\nn
\eeqa
are the Euler - Lagrange operators.

Recall that the densities of the local functionals are considered modulo total $x$-derivatives. So, the Poisson bracket \eqref{pb1} vanishes {\sl iff} the integrand in the rhs is a total $x$-derivative. The following important criterion (see, e.g., \cite{dt}) gives a necessary and sufficient condition for a differential polynomial
to be a total $x$-derivative of another differential polynomial.

\begin{theorem} Let
$$
h=\sum_{k,l\geq 0}^M\sum_{i, j} A_{i, j}(u,v)\, \pal_x^{i_1} u \dots \pal_x^{i_k} u \,\pal_x^{j_1} v\dots \pal_x^{j_l} v
$$
be a differential polynomial with coefficients $A_{i,j}(u,v)$ smooth on a disk ${\mathcal D}\subset \mathbb R^2$. Here
$$
i=(i_1, \dots, i_k), \quad j=(j_1, \dots, j_l), \quad 1\leq i_1, \dots, i_k, \, j_1, \dots, j_l \leq N
$$
are multiindices  of the summation, $M$, $N$ are two nonnegative integers. Then there exists another differential polynomial $\tilde h$ of the same form such that
$$
h=\pal_x \tilde h, \quad \pal_x =\sum_{m\geq 0} \left( u^{(m+1)} \frac{\pal}{\pal u^{(m)}} +  v^{(m+1)} \frac{\pal}{\pal v^{(m)}} \right)
$$
{\rm iff}
$$
\frac{\delta H}{\delta u(x)}=0, \quad \frac{\delta H}{\delta v(x)}=0,
\quad H=\int h\, dx.
$$
\end{theorem}

\begin{cor} Two local functionals $F=\int f\, dx$ and $G=\int g\, dx$ commute with respect to the Poisson bracket \eqref{pb} {\rm iff}
\eqa\label{kommut}
&&
{\rm E}_u \left( \frac{\delta F}{\delta u(x)} \,\pal_x 
\frac{\delta G}{\delta v(x)}+\frac{\delta F}{\delta v(x)}\, \pal_x 
\frac{\delta G}{\delta u(x)}\right)=0
\nn\\
&&
\\
&&
{\rm E}_v \left( \frac{\delta F}{\delta u(x)} \,\pal_x 
\frac{\delta G}{\delta v(x)}+\frac{\delta F}{\delta v(x)}\, \pal_x 
\frac{\delta G}{\delta u(x)}\right)=0.
\nn
\eeqa
\end{cor}

Alternatively \eqref{wave1} can be recast into the form of Euler - Lagrange equations
$$
\delta S=0, \quad S=\int\int \left[\frac12 \phi_t^2 - V(\phi_x)\right]\, dx\, dt
$$
with
$$
P(u)=V'(u)
$$
by introducing a potential
$$
\phi_x=u.
$$
The canonical Poisson bracket
$$
\{ \phi(x), \phi_t(y)\}=\delta(x-y)
$$
yields \eqref{pb}.

\subsection{Domain of hyperbolicity and Riemann invariants}\par

From \eqref{wave2} it readily follows that the nonlinear wave equation is hyperbolic on the domain of convexity of $P(u)$,
\beq\label{hyper1}
( u, v) \in \mathbb R^2 \quad \mbox{such that}\quad P''(u) >0
\eeq
and elliptic when $P(u)$ becomes concave. On the domain of hyperbolicity \eqref{wave2} can be reduced to the diagonal form
\beq\label{wave31}
\pal_t r_{\pm} =\pm \sqrt{P''(u)} \, \pal_x r_{\pm}
\eeq
introducing Riemann invariants
\beq\label{riem2}
r_{\pm}=v\pm Q(u), \quad \mbox{where}\quad Q'(u) =\sqrt{P''(u)}.
\eeq
We also remind the following property of solutions written in the Riemann invariants: introducing the {\it characteristic variables}
\beq\label{khar1}
x_\pm =x\pm \sqrt{P''(u)} t
\eeq
the function $r_+ / r_-$ does not depend on $x_-/ x_+$, 
\beq\label{khar2}
\frac{\pal r_+}{\pal x_-} =\frac{\pal r_-}{\pal x_+}=0.
\eeq

\subsection{Solutions of nonlinear wave equations}\par

Let us first describe the local structure of general solution to the system \eqref{wave2}.

\begin{lemma} \label{l-hodo} Let $f=f(u,v)$ be an arbitrary solution to the linear PDE
\beq\label{chap1}
f_{uu}=P''(u) f_{vv}.
\eeq
Let $(u_0, v_0)$ be a point such that 
\beq\label{jac1}
\left[f_{uv}(u_0, v_0)\right]^2-P''(u_0) \left[ f_{vv}(u_0, v_0)\right]^2\neq 0.
\eeq
Then the system of equations
\eqa\label{hodo1}
&&
x=f_u(u,v)
\nn\\
&&
\\
&&
t=f_v(u,v)
\nn
\eeqa
defines a unique pair of functions $\left( u(x,t), v(x,t)\right)$ near the point
$$
\left( x_0, t_0\right):= \left( f_u(u_0, v_0), f_v(u_0, v_0)\right)
$$
such that
$$
u(x_0, t_0)=u_0, \quad v(x_0, t_0)=v_0.
$$
The functions $\left( u(x,t), v(x,t)\right)$ solve \eqref{wave2}.
Their $x$-derivatives at the point $(x_0, t_0)$ do not vanish simultaneously. Conversely, any solution to \eqref{wave2}
satisfying the condition
\beq\label{monoton}
\left[v_x^2 -P''(u) u_x^2\right]_{(x_0, t_0)}\neq 0
\eeq
can be locally obtained by the above procedure starting from a suitable solution $f(u,v)$ to the linear PDE \eqref{chap1}.
\end{lemma}

This is a reformulation of the classical method of characteristics (also called hodograph transformation).

\begin{remark} The Riemann invariants \eqref{riem2} are also the characteristic variables for the linear PDE \eqref{chap1}, i.e., in these variables the PDE for the function $f=f(r_+, r_-)$ reads
\beq\label{petrov}
\frac{\pal^2 f}{\pal r_+ \pal r_-} -\frac18\, \frac{P'''}{(P'')^{3/2}} \left( \frac{\pal f}{\pal r_+} -\frac{\pal f}{\pal r_-}\right)=0.
\eeq
\end{remark}

\begin{remark} The exceptional solutions to \eqref{wave2}
that cannot be written in the form \eqref{hodo1} can be represented in the form
$$
v=Q(u) + c \quad \mbox{\rm or}\quad v=-Q(u)+c
$$
for a constant $c$ where $u=u(x,t)$ satisfies the {\rm scalar} quasilinear PDE
$$
u_t =\sqrt{P''(u)}\, u_x \quad \mbox{\rm or}\quad  u_t =-\sqrt{P''(u)} \,u_x.
$$
On these solutions one of the Riemann invariants takes constant values.
\end{remark}

\subsection{Conserved quantities and infinitesimal symmetries
of nonlinear wave equation}\par

An alternative interpretation of solutions to the linear
PDE \eqref{chap1} is given by

\begin{lemma} Given a function $f=f(u,v)$, consider the Hamiltonian
\beq\label{hf0}
H_f=\int f(u,v)\, dx.
\eeq
This Hamiltonian commutes with the Hamiltonian \eqref{hwave} of nonlinear wave equation {\rm iff} the function $f$ satisfies \eqref{chap1}.
\end{lemma}

{\sl Proof}. The Poisson bracket \eqref{pb} of $H_f$ and $H$ reads
$$
\{ H_f, H\} =\int \left[ \frac{\delta H_f}{\delta u(x)} \,\pal_x 
\frac{\delta H}{\delta v(x)}+\frac{\delta H_f}{\delta v(x)}\, \pal_x 
\frac{\delta H}{\delta u(x)}\right]\, dx
=\int\left[ f_u v_x + f_v P''(u) u_x\right]\, dx.
$$
The integrand is a total $x$-derivative {\sl iff}
$$
\pal_u \left( f_u\right) =\pal_v\left(f_v P''(u)\right).
$$
This gives \eqref{chap1}. \epf

Recall that the differential polynomial $f=f(u,v; u_x, v_x, \dots)$
is called {\it conserved quantity} (or density of a {\it first integral})  of the PDE \eqref{wave2} if the following identity holds
\beq\label{conserv2}
\pal_t f =\pal_x g, \quad \pal_t= \sum_{m\geq 0} \left[ \left( u_t\right)^{(m)} \frac{\pal}{\pal u^{(m)}}+ \left( v_t\right)^{(m)} \frac{\pal}{\pal v^{(m)}}\right]
\eeq
for some differential polynomial $g$. It is understood that in \eqref{conserv2} one has to substitute from \eqref{wave2}
$u_t\mapsto v_x$, $v_t \mapsto P''(u) u_x$.

\begin{cor} The function $f=f(u,v)$ is a conserved quantity for the nonlinear wave equation \eqref{wave2} {\rm iff} it satisfies the linear PDE \eqref{chap1} .
\end{cor}

{\sl Proof}. For a solution $f$ to \eqref{chap1} one has
$$
\pal_t f =\pal_x g
$$
where $g=g(u,v)$ is determined by a quadrature
$$
dg= P''(u) f_v du + f_u dv.
$$
The converse statement is straightforward. \epf

\begin{cor} The Hamiltonian PDE
\eqa\label{symwave}
&&
u_s=\pal_x f_v(u,v)
\nn\\
&&
\\
&&
v_s =\pal_x f_u(u,v)
\nn
\eeqa
commutes with the wave equation
$$
\left( u_t\right)_s=\left( u_s\right)_t, \quad \left( v_t\right)_s=\left( v_s\right)_t
$$
if the function $f(u,v)$ satisfies the linear PDE \eqref{chap1}.
\end{cor}

{\sl Proof}. Indeed, the flow \eqref{symwave} is generated by the Hamiltonian $H_f$:
\eqa
&&
u_s=\pal_x \frac{\delta H_f}{\delta v(x)}
\nn\\
&&
v_s=\pal_x \frac{\delta H_f}{\delta u(x)}.
\nn
\eeqa
As it is well known from Hamiltonian mechanics, commutativity of Hamiltonians implies commutativity of the flows. \epf

\subsection{Commutative Lie algebra of infinitesimal symmetries of  nonlinear wave equation
}\par
We conclude that the solutions to the linear PDE \eqref{chap1} are associated with infinitesimal symmetries of the nonlinear wave equation \eqref{wave1}. 
A somewhat stronger statement says that all Hamiltonians of the form $H_f$ commute pairwise, i.e., the Lie algebra of symmetries is commutative.

\begin{lemma} Given two arbitrary solutions $f$, $g$ to the linear PDE \eqref{chap1}, the Hamiltonians
$$
H_f=\int f(u,v)\, dx \quad \mbox{\rm and}\quad H_g=\int g(u,v)\, dx
$$
commute:
$$
\{ H_f, H_g\}=0.
$$
\end{lemma}

We leave the proof of this lemma as an exercise for the reader.

For a given $P(u)$ one obtains therefore an infinite family of commuting Hamiltonians parametrized by solutions to the linear PDE \eqref{chap1}. The associated Hamiltonian systems commute pairwise.

\begin{remark} All equations \eqref{symwave} of the commutative family have the same Riemann invariants \eqref{riem2}:
\beq\label{riem3}
\pal_s r_{\pm} =(h_{uv}\pm \sqrt{P''(u)} h_{vv}) \,\pal_x r_{\pm}.
\eeq
\end{remark}

\begin{exam} Given a real number $\kappa\neq 0, \, \pm 1$, choose
\beq\label{poly1}
P(u) =\frac{u^\kappa}{\kappa (\kappa-1)}.
\eeq
The function
\beq\label{poly2}
f=\frac12 u\, v^2 +\frac{u^\kappa}{\kappa (\kappa+1)}
\eeq
satisfies the linear equation \eqref{chap1}. The associated Hamiltonian flow reads
\eqa\label{poly3}
&&
u_t =\pal_x (u\, v)
\nn\\
&&
\\
&&
v_t =\pal_x\left( \frac{v^2}2 + \frac{u^\kappa}{\kappa}\right).
\nn
\eeqa
These equations coincide (after the change of sign of the time variable) with equations of motion of one-dimensional isentropic gas with the equation of state
$$
p\sim \rho^{\kappa+1}.
$$
Here $u=\rho$ (the mass density), $v$ the velocity. Therefore the equations \eqref{poly3} of gas dynamic can be considered
as an infinitesimal symmetry of the nonlinear wave equation
$$
u_{tt}=\frac{u^{\kappa-1}}{\kappa-1} u_{xx} + u^{\kappa-2} u_x^2.
$$
\end{exam}

\begin{exam} For the exceptional value $\kappa=1$ take
\beq\label{nls1}
P(u) =\mp u\, (\log u -1).
\eeq
The choice
\beq\label{nls2}
f=\frac12 \left( u^2 \mp u\, v^2\right)
\eeq
of a solution to the linear PDE \eqref{chap1} produces the dispersionless limit 
\eqa\label{nlsd}
&&
u_t +(u\, v)_x=0
\\
&&
v_t +v\, v_x \mp u_x =0
\nn
\eeqa
of the focusing/defocusing nonlinear Schr\"odinger equation
\beq\label{nls0}
i\, \psi_t +\frac12\, \psi_{xx} \pm |\psi|^2\psi=0
\eeq
written in the coordinates
$$
u=|\psi|^2, \quad v=\frac1{2i}\left( \frac{\psi_x}{\psi} -\frac{\bar\psi_x}{\bar\psi}\right),
$$
i.e.,
\eqa\label{nls}
&&
u_t +(u\, v)_x=0
\nn\\
&&
v_t +v\, v_x \mp u_x =\frac14 \left( \frac{u_{xx}}{u}-\frac12\, \frac{u_x^2}{u^2}\right)_x.
\nn
\eeqa
So, the dispersionless NLS is an infinitesimal symmetry of the following nonlinear wave equation
$$
u_{tt}\pm \pal_x^2 \left( \frac1{u}\right)=0.
$$
\end{exam}

\begin{remark}\label{rem_fro} In \cite{npb} a hierarchy of the first order commuting quasilinear PDEs has been associated with an arbitrary solution $F=F(v^1, \dots, v^n)$ of equations of associativity
\eqa\label{wdvv}
&&
\frac{\pal^3 F}{\pal v^i \pal v^j \pal v^p}\eta^{pq} \frac{\pal^3 F}{\pal v^q \pal v^k\pal v^l}=\frac{\pal^3 F}{\pal v^l \pal v^j \pal v^p}\eta^{pq} \frac{\pal^3 F}{\pal v^q \pal v^k\pal v^i}, \quad i, 
j, k, l=1, \dots, n
\nn\\
&&
\\
&&
\frac{\pal^3 F}{\pal v^1 \pal v^i \pal v^j}=\eta_{ij}.
\nn
\eeqa
Here $\eta^{ij}$ are entries of a constant symmetric nondegenerate matrix, $\eta_{ij}$ are entries of the inverse matrix. The quasihomogeneity axiom has not been used for the construction of the hierarchy. For $n=2$ the solution $F=F(v,u)$ must have the following form
$$
F=\frac12 \, v\, u^2 + \Phi(u)
$$
for some function $\Phi(u)$. For
$$
\Phi'(u)=P(u)
$$
this family of commuting PDEs coincides with the one associated with solutions to \eqref{chap1}. The Hamiltonians of the hierarchy of \cite{npb} are solutions to \eqref{chap1} polynomial in $v$.
\end{remark}

Before proceeding to the discussion of the critical phenomena let us extend the method of characteristics described in Lemma \ref{l-hodo} to solving any equation of the commutative family
\eqa\label{h-eq}
&&
u_s=\pal_x h_v
\nn\\
&&
\\
&&
v_s =\pal_x h_u
\nn
\eeqa
with the Hamiltonian
$$
H_h =\int h(u,v)\, dx, \quad h_{uu}=P''(u) h_{vv}.
$$

\begin{lemma} Let $f=f(u,v)$ be any solution to \eqref{chap1} satisfying the condition
\beq
\left( f_{v}^0 h_{uvv}^0 - f_{uv}^0 h_{vv}^0\right)^2 -P''(u_0) 
\left( f_{v}^0 h_{vvv}^0 - f_{vv}^0 h_{vv}^0\right)^2\neq 0
\nn
\eeq
for some $(u_0, v_0)$.
 Then the system
\eqa\label{hodo2}
&&
x+s\, h_{uv}(u,v) = f_u(u,v)
\nn\\
&&
\\
&&
\quad\quad s\, h_{vv}(u,v) = f_v(u,v)
\nn\\
&&
\nn\\
&&
u(x_0, s_0)=u_0, \quad v(x_0, s_0)=v_0
\nn
\eeqa
defines a unique pair of functions $(u(x,s), v(x,s))$ near $(x_0, s_0)$ determined from
\eqa
&&
x_0+s_0\, h_{uv}(u_0,v_0) = f_u(u_0,v_0)
\nn\\
&&
\quad\quad  s_0\, h_{vv}(u_0,v_0) = f_v(u_0,v_0).
\nn
\eeqa
These functions solve the system \eqref{h-eq}. Conversely, any solution to \eqref{h-eq} satisfying certain nondegeneracy conditions can be obtained by this procedure.
\end{lemma}

\setcounter{equation}{0}
\setcounter{theorem}{0}
\section{Singularities of solutions to nonlinear wave equation}\label{sect3}\par

In this section we develop a local classification  of singularities of solutions to the nonlinear wave equation. Actually, the technique based on the implicit function theorem works equally well also for any PDE of the form \eqref{h-eq} from the commuting family of infinitesimal symmetries of the wave equation. So, our study of singularities will be developed for the solutions to the general system \eqref{h-eq}.

The functions $\left( u(x,s), v(x,s)\right)$ determined by the system \eqref{hodo2} are smooth provided the conditions of the implicit function theorem hold true. At the points $(x_{\rm c}, s_{\rm c})$ where the implicit function theorem fails to be applicable the solution 
to \eqref{h-eq} has a {\it weak singularity}, i.e., there exist the limits
$$
\lim_{x\to x_{\rm c}, \, s\to s_{\rm c}} u(x,s)=u_{\rm c}, \quad \lim_{x\to x_{\rm c}, \, s\to s_{\rm c}} v(x,s)=v_{\rm c}
$$
but the derivatives $u_x$, $v_x$ blow up. The main goal of this Section is to describe the local behaviour of a solution near the {\it generic} weak singularity. The precise definition of genericity will be given below.

Let us first subdivide weak singularities into three classes.

Type I. The limiting value $(u_{\rm c}, v_{\rm c})$ belongs to the domain of hyperbolicity of the system \eqref{h-eq}, i.e.
\beq
\label{typ1}
P''(u_{\rm c})>0.
\eeq

Type II. The limiting value $(u_{\rm c}, v_{\rm c})$ belongs to the domain of ellipticity of the system \eqref{h-eq}, i.e.
\beq
\label{typ2}
P''(u_{\rm c})<0.
\eeq

Type III. The limiting value $(u_{\rm c}, v_{\rm c})$ belongs to the boundary between the domains of hyperbolicity and ellipticity,
\beq\label{typ3}
P''(u_{\rm c})=0.
\eeq

Observe that the linear PDE \eqref{chap1} belongs to the hyperbolic/elliptic type near the critical point $(u_{\rm c}, v_{\rm c})$
of Type I/II, while the type of this PDE changes near a generic (i.e., $P'''(u_{\rm c})\neq 0$) critical point of Type III. 

In this Section we will study the local behaviour of the solution near the weak singularity following the above classification. Only Type I and Type II singularities will be considered; the study of Type III singularities is postponed to a subsequent publication.
 
Let us do first some general simplification of notations. First, since the partial derivative $h_v$ satisfies the linear PDE \eqref{chap1} if the function $h$ does so, it suffices to study singularities of implicit functions $u(x,s)$ and $v(x,s)$ determined by the system
\eqa\label{imp}
&&
x+s\, h_u(u,v) = f_u(u,v)
\nn\\
&&
\\
&&
\quad\quad s\, h_v(u,v)= f_v(u,v)
\nn
\eeqa
where both functions $h$ and $f$ satisfy \eqref{chap1}. Next, the shifts
$$
x\mapsto x-x_{\rm c}, \quad s\mapsto s-s_{\rm c}
$$
can be absorbed by a redefinition
$$
f\mapsto f+ x_{\rm c} u + s_{\rm c} h
$$
of the function $f$. So, without loss of generality we may assume that $x_{\rm c}=s_{\rm c}=0$.

We will also adopt the following system of notations: the values
of the functions $f(u,v)$, $h(u,v)$  and of their derivatives $f_u(u,v)$, $f_v(u,v)$ etc. at the critical point $(u_{\rm c}, v_{\rm c})$ will be denoted $f^0$, $h^0$, $f^0_u$, $f_v^0$ etc. We will also denote
$$
P_0:= P(u_{\rm c}), \quad P'_0 :=P'(u_{\rm c}) 
$$
etc.

According to the above conventions we will have
\beq\label{crit1}
f_u^0=f_v^0=0.
\eeq
Violation of the implicit function theorem conditions in our notations reads
\beq\label{crit2}
-\det\left( \begin{array}{cc} f_{uu}-s \, h_{uu} & f_{uv}-s \, h_{uv}\\
f_{uv}-s \, h_{uv} & f_{vv}-s \, h_{vv}\end{array}\right)_{x=0, \, s=0} = \left( f_{uv}^0\right)^2 -P''_0 \left( f_{vv}^0\right)^2=0.
\eeq

\subsection{Critical points of Type I and Whitney singularities}\par 

As it follows from the representation \eqref{hodo1} of solutions, the inverse map to
$$
(x,t)\mapsto \left(u(x,t), v(x,t)\right)
$$
is well defined also on a neighborhood of the critical point. Moreover, it is given by the gradient
\beq\label{grad}
(u,v)\mapsto \left( f_u(u,v), f_v(u,v)\right)=(x,t)
\eeq
of a function satisfying the linear PDE \eqref{chap1} that, in the case of Type I singularity at $(u_{\rm c}, v_{\rm c}$) belongs to the hyperbolic type. The points of blow-up of the derivatives of the solution $\left(u(x,t),v(x,t)\right)$ are exactly the points where the Jacobian of the gradient map \eqref{grad} vanishes.

Classification of singularities of smooth maps of plane to plane 
obtained in the celebrated paper of H.Whitney \cite{w} was one of the first successes of the singularity theory. Recall that Whitney obtained the complete list of normal forms of smooth maps
\eqa
&&
(\mathbb R, 0)\to (\mathbb R, 0)
\nn\\
&&
(u,v) \mapsto (x,t)
\nn
\eeqa
stable with respect to small perturbations. He proved that, up to diffeomorphisms
\eqa
&&
(x,t)\mapsto (\tilde x, \tilde t)
\nn\\
&&
(u,v)\mapsto (\tilde u, \tilde v)
\nn
\eeqa
 these maps belong to one of the following three normal forms:

\noindent Type $W1$
\eqa\label{whit1}
&&
\tilde x= \tilde v
\nn\\
&&
\\
&&
\tilde t=\tilde u
\nn
\eeqa

\noindent Type $W2$
\eqa\label{whit2}
&&
\tilde x = \tilde v^2
\nn\\
&&
\\
&&
\tilde t= \tilde u
\nn
\eeqa

\noindent Type $W3$
\eqa\label{whit3}
&&
\tilde x=\tilde u\, \tilde v-\frac16 \tilde v^3
\nn\\
&&
\\
&&
\tilde t=\tilde u
\nn
\eeqa

For the $W1$ case the inverse functions $u=u(x,t)$, $v=v(x,t)$ are well defined and smooth.
We will show below that also in the general situation the solutions to the quasilinear system \eqref{h-eq} near a generic Type I singularity locally has the universal shape \eqref{whit2}. 

Another natural question is about the structure of solutions near {\it the first} point of singularity.
We say that $(x_{\rm c}, t_{\rm c})$ is the first singularity of the solution to the nonlinear wave equation if the functions $u(x,t)$ and $v(x,t)$ are smooth
for sufficiently small $|x-x_{\rm c}|$ and for $t<t_{\rm c}$. We will see that the local behaviour of a generic solution near a first point singularity, within the hyperbolicity domain, is described by the $W3$ singularity on the Whitney list. This sounds like being almost obvious; an additional point is that only {\it affine transformation} of independent variables $x$, $t$ will be allowed near the critical point.

Let proceed now to the study of singularities of solutions to the general system \eqref{h-eq}. Due to the assumption $P''_0>0$ the jacobian \eqref{crit2} factorizes
$$
 \left( f_{uv}^0\right)^2 -P''_0 \left( f_{vv}^0\right)^2=\left( f_{uv}^0 + \sqrt{P''_0} f_{vv}^0\right)\left( f_{uv}^0 - \sqrt{P''_0} f_{vv}^0\right).
 $$
 The first of the genericity assumption says that, at the critical point only one of the factors vanishes. Without loss of generality we may assume that
\beq\label{crit11}
f_{uv}^0 - \sqrt{P''_0} f_{vv}^0=0.
\eeq

\begin{defi} We say that the critical point \eqref{crit11} of Type 1 is {\rm generic} if:
\beq\label{generic11}
f_{uv}^0 + \sqrt{P''_0} f_{vv}^0\neq 0
\eeq
and
\beq\label{generic12}
f_{uvv}^0 -\sqrt{P_0''} f_{vvv}^0 + f_{vv}^0\, \frac{P_0'''}{4\, P_0''}\neq 0
\eeq
or
$$
f_{uvv}^0 -\sqrt{P_0''} f_{vvv}^0 + f_{vv}^0\, \frac{P_0'''}{4\, P_0''}= 0 
$$
but
\eqa\label{generic13}
&& 
A_0:= \quad h_{uv}^0 -\sqrt{P_0''} h_{vv}^0 + h_v^0 \, \frac{P_0'''}{4\, P_0''}\neq 0 \quad \mbox{\rm and}
\nn\\
&&
 \\
&&
B_0:=f_{uvvv}^0 -\sqrt{P_0''}\, f_{vvvv}^0 + f_{vvv}^0 \, \frac{P_0'''}{4\, P_0''} +\frac1{32} \frac{f_{vv}^0}{\left( P_0''\right)^{5/2}} \left[ 5 (P_0''')^2 -4 P_0^{IV} P_0''\right]\neq 0.
\nn
\eeqa
\end{defi}

The formula
\beq\label{drx}
\pal_x r_{\pm} = \frac1{f_{uv}^0 \pm \sqrt{P_0''} f_{vv}^0 -s\,\left(
h_{uv}^0 \pm \sqrt{P_0''} h_{vv}^0\right)}
\eeq
for the $x$-derivatives of Riemann invariants clarifies the meaning of the first genericity assumption \eqref{generic11}: at the critical point $(x_{\rm c}, s_{\rm c})=(0,0)$ only one Riemann invariant blows up (under our assumptions $r_-$ does so) while another one remains smooth.

We now show that the case \eqref{generic12} is incompatible with the assumption that $(x_{\rm c}, s_{\rm c})=(0,0)$ is  locally the first critical point. 

\begin{lemma}\label{lm32} Let $(x_{\rm c}, s_{\rm c})=(0,0)$ be a critical point of Type I satisfying \eqref{crit11}, \eqref{generic11} such that the functions $u(x,s)$, $v(x,s)$ determined by \eqref{imp} along with normalization $u(x_{\rm c}, s_{\rm c})=u_{\rm c}$, $v(x_{\rm c}, s_{\rm c})=v_{\rm c}$ are smooth for $s<s_{\rm c}$ for sufficiently small $|x-x_{\rm c}|$, $|s-s_{\rm c}|$. Then
\beq\label{generic14}
f_{uvv}^0 -\sqrt{P_0''} f_{vvv}^0 + f_{vv}^0\, \frac{P_0'''}{4\, P_0''}= 0.
\eeq
\end{lemma}

The meaning of the constraint \eqref{generic14} is the following: the graph of the Riemann invariant $r_-(x, s_{\rm c})$ must have an inflection point at $x=x_{\rm c}$.

{\sl Proof}. For simplicity of calculations we will assume that
$$
h_v^0\neq 0.
$$
Under this assumption one can solve the system \eqref{imp} in the form
$$
s=s(u,v)=\frac{f_v}{h_v}, \quad x=x(u,v)=\frac{f_u h_v-f_v h_u}{h_v}.
$$
Expanding these functions in Taylor series  in $r_+$, $r_-$ at $u=u_{\rm c}$, $v=v_{\rm c}$ we obtain the following expressions for the linear combinations\footnote{Clearly these correspond to the characteristics of the hyperbolic PDE.}
\beq\label{xpm}
x_{\pm}=x+\left( h_u^0 \pm \sqrt{P_0''} h_v^0\right)\, s
\eeq
of independent variables
\eqa\label{riad12}
&&
x_+ = 2 f_{vv}^0\sqrt{P_0''}\, \bar r_+ 
\nn\\
&&
\qquad+\frac1{2\, h_v^0} \left[ h_v^0 \left( f_{uvv}^0 + \sqrt{P_0''}\, f_{vvv}^0\right) -2 f_{vv}^0 \left( h_{uv}^0 + \sqrt{P_0''}\, h_{vv}^0\right)\right] \bar r_+^2 +{\mathcal O}(|w|^3)
\nn\\
&&
\\
&&
x_- = f_{vv}^0 \,\frac{P_0'''}{8 P_0''} \, \bar r_+^2 
-\frac{f_{vv}^0}{h_v^0} \left( h_{uv}^0 -\sqrt{P_0''}\, h_{vv}^0 + h_v^0 \,\frac{P_0'''}{4 P_0''} \right) \bar r_+ \bar r_-
\nn\\
&&
\qquad +\frac12 \left( f_{uvv}^0 -\sqrt{P_0''} f_{vvv}^0 + f_{vv}^0\, \frac{P_0'''}{4\, P_0''}\right) \bar r_-^2 +{\mathcal O}(|w|^3)
\nn
\eeqa
where we introduced the notations
\beq\label{wpm}
\bar r_{\pm} = r_\pm -r_\pm^0
\eeq
for the shifted Riemann invariants.
Assuming \eqref{generic12} do a rescaling
$$
x_{\pm} \mapsto k\, x_{\pm}, \quad \bar r_+ \mapsto k\, \bar r_+, \quad \bar r_- \mapsto k^{1/2} \bar r_-
$$
with a small positive $k$. After simple calculations one obtains
\eqa\label{whi11}
&&
\tilde x_+ =  \bar r_+ +{\mathcal O}(k^{1/2})
\nn\\
&&
\\
&&
\tilde x_- =  \bar r_-^2 +{\mathcal O}(k^{1/2})
\nn
\eeqa
where
\eqa
&&
\tilde x_+   =\frac{x_+}{2 f_{vv}^0\sqrt{P_0''}} 
\nn\\
&&
\nn\\
&&
\tilde x_- = \frac{2 \, x_-}{f_{uvv}^0 -\sqrt{P_0''} f_{vvv}^0 + f_{vv}^0\, \frac{P_0'''}{4\, P_0''}}
\nn
\eeqa
(observe that $f_{vv}^0\neq 0$ due to the genericity assumption \eqref{generic11}).
We obtain Whitney normal form of type $W2$. Therefore the solution to the implicit function system \eqref{imp} near the critical point is defined on a half plane with the boundary $\tilde x_-=0$ passing through the critical point. Clearly such a behaviour is incompatible with the assumptions of the lemma. \epf

\begin{theorem} Under the genericity assumption \eqref{generic13} the local structure of the solution to the implicit function equations is described by the Whitney normal form of type $W3$.
\end{theorem}

\pf We have to include the degree 3 terms in the Taylor
expansions of $x(u,v)$, $s(u,v)$. The rescaling
\beq\label{scal2}
x_+\mapsto k^{2/3} x_+, \quad x_- \mapsto k\, x_-, \quad \bar r_+ \mapsto k^{2/3} \bar r_+, \quad \bar r_- \mapsto k^{1/3} \bar r_-
\eeq
yields
\eqa\label{tip1}
&&
x_+ = 2 f_{vv}^0 \sqrt{P_0''}\, \bar r_+ +{\mathcal O}(k^{1/3})
\nn\\
&&
\\
&&
x_-= - \frac{f_{vv}^0}{h_v^0} A _0\, \bar r_+ \bar r_- +\frac16 B_0\, \bar r_-^3 
+{\mathcal O}(k^{1/3})
\nn
\eeqa
where $A_0$ and $B_0$ are defined in \eqref{generic13}. The theorem is proved. \epf

\subsection{Critical points of Type II and elliptic umbilic catastrophe}

Recall \cite{th} that elliptic umbilic catastrophe can be represented as the singularity at $X=Y=0$ of the function
\beq\label{umb}
F(X,Y)=\frac13 \,X^3 - X\, Y^2.
\eeq
It can also be identified (see \cite{ar}) as one of the real forms of the codimension 2 singularity of type $D_4$. The following interpretation will be useful in sequel. Let $\Phi(X,Y)$ be a harmonic function of two variables,
$$
\Phi_{XX}+\Phi_{YY}=0.
$$
Consider singularities of the gradient map
\beq\label{grad1}
(X,Y)\mapsto \left( \Phi_X(X,Y), \Phi_Y(X,Y)\right) =(U,V)
\eeq
i.e. the points where the Jacobian of \eqref{grad1} vanishes,
$$
\det \left( \begin{array}{cc} \pal U/\pal X & \pal U\pal Y\\
 & \\
\pal V/\pal X & \pal V/\pal Y\end{array}\right) =-\left[ \Phi_{XX}^2 +\Phi_{XY}^2\right]=0.
$$
As it follows from the above formula for the Jacobian, at this point all the second derivatives of $\Phi$ vanish. So, the 3-jet of the function must be a cubic harmonic polynomial. It must coincide with \eqref{umb} up to a linear change of variables.

The gradient map of \eqref{umb} is given by the following system
\beq\label{umb1}
\left.\begin{array}{cc} U & =X^2 -Y^2\\
 & \\
V & =-2\, X\, Y\end{array}\right\} .
\eeq
In the complex coordinates
$$
W=U-i\, V, \quad Z=X+i\, Y
$$
this gives the quadratic map
\beq\label{umb2}
W= Z^2.
\eeq

We will now show that the Type II singularity of a generic solution to the quasilinear system \eqref{h-eq} is locally described by the elliptic umbilic catastrophe.

As above we assume that the singular point is $x_C=s_C=0$; it will also be assumed, as above, that
$$
f_u^0=f_v^0=0.
$$ 
Denote $\bar u =u- u_C$, $v=v-v_C$. Consider a solution to \eqref{h-eq} determined by the system \eqref{hodo2} near the critical point.

\begin{lemma} At the critical point all second derivatives
of $f$ vanish:
\beq\label{d2}
f_{uu}^0=f_{vv}^0=f_{uv}^0=0.
\eeq
\end{lemma}

{\sl Proof} follows from \eqref{crit2} due to negativity of $P_0''$.

\begin{defi} The Type II critical points is called {\rm generic} if the 3-jet of $f(u,v)$ does not vanish at this point.
\end{defi}

\begin{lemma} After the rescaling
$$
\begin{array}{clr} x & \mapsto  & k\, x\\
s & \mapsto & k\, s\\
\bar u & \mapsto  & k^{1/2} \bar u\\
\bar v & \mapsto & k^{1/2} \bar v\end{array}
$$
at the limit $k\to 0+$ the solution to the system \eqref{crit2}
near a generic Type II critical point $(x=x_C=0, s=s_C=0, u=u_C, v=v_C)$ tends to the solution of the equation
\eqa\label{umb3}
&&
z=\frac12\, a_0 w^2 
\\
&&
\nn\\
&&
z:=x+s\, h_u^0 +i\, c_0 s\, h_v^0, \quad w:=\bar v +i\, c_0\, \bar u , \quad a_0:= f_{uvv}^0 +i\, c_0\, f_{vvv}^0
\nn
\eeqa
where
$$
c_0=\sqrt{-P_0''}.
$$
\end{lemma}

Observe that, due to \eqref{chap1} and \eqref{d2} at the critical point the third derivatives of the function $f(u,v)$ satisfy
$$
f_{uuu}^0=P_0'' f_{uvv}^0, \quad f_{uuv}^0=P_0'' f_{vvv}^0.
$$
So, at the generic Type II critical point we have necessarily $a_0\neq 0$.

{\sl Proof} of Lemma is similar to that of Lemma \ref{lm32}. Expanding the left hand sides of eqs. \eqref{imp} in Taylor series and using equation \eqref{chap1} valid for both functions $f$ and $h$  one obtains
\eqa
&&
x+s\, h_u^0 +s\, (-c_0^2 h_{vv}^0 \bar u + h_{uv}^0 \bar v) +\frac{s}2 \left[ h_{uvv}^0 (-c_0^2 \bar u^2 +\bar v^2) -2 c_0^2 h_{vvv}^0\bar u\, \bar v\right]
\nn\\
&&
\nn\\
&&=
\frac12\left[ f_{uvv}^0 (-c_0^2 \bar u^2 + \bar v^2) -2 c_0^2 f_{vvv}^0 \bar u\, \bar v\right] + {\mathcal O}\left((\bar u^2 +\bar v^2)^{3/2}\right)
\nn\\
&&
\nn\\
&&
s\, h_v^0 + s\, (h_{uv}^0 \bar u + h_{vv}^0 \bar v) +\frac{s}2 \left[ h_{vvv}^0 (-c_0^2 \bar u^2 +\bar v^2) +2 h_{uvv}^0 \bar u\, \bar v\right]
\nn\\
&&
\nn\\
&&
=\frac12 \left[ f_{vvv}^0 (-c_0^2 \bar u^2 + \bar v^2) +2f_{uvv}^0 \bar u\, \bar v\right] + {\mathcal O}\left((\bar u^2 +\bar v^2)^{3/2}\right)
\nn
\eeqa
After the rescaling and division by $k$ one arrives at the system
$$
\begin{array}{rc} x +s\, h_{u}^0 & = \frac12 \left[ f_{uvv}^0 (-c_0^2 \bar u^2 + \bar v^2)-2 c_0^2 f_{vvv}^0 \bar u\, \bar v\right] +{\mathcal O}\left( k^{1/2}\right)\\
 & \\
 s\, h_v^0 & = \frac12 \left[ f_{vvv}^0 (-c_0^2 \bar u^2 + \bar v^2) + 2 f_{uvv}^0 \bar u\, \bar v\right]  +{\mathcal O}\left( k^{1/2}\right).\end{array}
 $$
This converges to the quadratic equation \eqref{umb3}. The Lemma is proved.

\medskip 

We see that, after a linear change of independent variables the local structure of the solution near a generic Type II singularity is described by the bifurcation diagram \eqref{umb2} of the elliptic umbilic catastrophe.

\subsection{Generalization to the multicomponent case}\par
Let us consider an $n$-component system of PDEs written in the diagonal form
\beq\label{semiham1}
\pal_t u_i=a_i(\ub)\pal_xu_i, \quad i=1, \dots, n,
\eeq
$n\geq 3$. Thus the dependent variables $u_1$, \dots, $u_n$ are the Riemann invariants for the system \eqref{semiham1}; the characteristic velocities $a_1(\ub)$, \dots, $a_n(\ub)$ will be assumed to be smooth functions taking {\it pairwise distinct} values,
$$
a_i(\ub)\neq a_j(\ub), \quad i\neq j\quad
\mbox{for all}\quad \ub \in {\mathcal D}.
$$
Recall the following definition \cite{tsarev}.

\begin{defi} The system \eqref{semiham1} is called {\rm semihamiltonian} if the coefficients satisfy the following differential equations
\beq\label{semiham2}
\pal_k\left( \frac{a_{i,j}}{a_i-a_j}\right) = \pal_j\left( \frac{a_{i,k}}{a_i-a_k}\right)
\eeq
for all triples of pairwise distinct indices $i$, $j$, $k$.
\end{defi}

Here
$$
\pal_i:=\frac{\pal}{\pal u_i}, \quad a_{i,j}:= \pal_j a_i;
$$
we will also use similar notations for higher derivatives, i.e.
$$
a_{i, jk}:= \pal_j
\pal_k a_i, \quad A_{i,jkl}:= \pal_j\pal_k\pal_l A_i
$$
etc. 

The meaning of this definition is clarified by the following theorem \cite{tsarev}.

\begin{theorem} {\rm (i)} Let $A_1(\ub)$, \dots, $A_n(\ub)$ be an arbitrary solution to the overdetermined system of PDEs
\beq\label{tsar1}
A_{i,j}= \frac{a_{i,j}}{a_i-a_j}\, \left(A_i -A_j\right), \quad i\neq j.
\eeq
Then the system of PDEs
\beq\label{semiham3}
\pal_su_i = A_i(\ub) \pal_x u_i, \quad i=1, \dots, n
\eeq
commutes with \eqref{semiham1},
\beq\label{comm}
\pal_s\pal_tu_i \equiv \pal_t\pal_su_i.
\eeq
{\rm (ii)} For a semihamiltonian system \eqref{semiham1} the general solution to \eqref{semiham3} depends on $n$ arbitrary functions of one variable.
\newline\noindent {\rm (iii)} Given a solution to the system \eqref{tsar1} such that
$$
A_1(\ub_0)=\dots = A_n(\ub_0)=0, \quad\prod_{i=1}^n A_{i,i}(\ub_0)\neq 0
$$
for some $\ub_0\in{\mathcal D}$, then the 
the system of equations 
\beq\label{tsar2}
x=a_i(\ub)\, t +A_i(\ub), \quad i=1, \dots, n
\eeq
defines a unique smooth vector function $\ub(x,t)$ for sufficiently small
$|x|$, $|t|$ such that $\ub(0,0)=\ub_0$. This function satisfies the system \eqref{semiham1}. Conversely, any solution $\ub(x,t)$ to the semihamiltonian system \eqref{semiham1} satisfying
$$
\pal_x u_1(0,0)\neq 0, \dots, \pal_x u_n(0,0)\neq 0
$$
can be locally obtained by the above procedure.
\end{theorem}

\begin{remark} For a semihamiltonian system \eqref{semiham1}
there also exists an infinite family of conserved quantities $h(\ub)$ such that the following identity holds true on the solutions to \eqref{semiham1}
\beq\label{conserv} 
\pal_t h(\ub) =\pal_x f(\ub)
\eeq
for some function $f(\ub)$. The conserved quantities are determined from the following overdetermined system of PDEs
\beq\label{tsar3}
\pal_i\pal_j h=\frac{a_{i,j}}{a_j-a_i}\,\pal_i h + \frac{a_{j,i}}{a_i-a_j}\, \pal_jh, \quad i\neq j.
\eeq
The solutions to this system also depend on $n$ arbitrary functions of one variable [{\it ibid}]. The conserved quantities $h(\ub)$ were called {\rm entropies} in \cite{serre}; the systems \eqref{semiham1} possessing a rich families of entropies (i.e., families depending on $n$ arbitrary functions of one variable) were called {\rm syst\`emes rich} in \cite{serre}. This class coincides with the class of semihamiltonian systems.
\end{remark}

Let us describe the local structure of singularities of solutions to the semihamiltonian systems near the point of gradient catastrophe. Without loss of generality we may assume that the gradient catastrophe takes place at the point $x_0=0$, $t_0=0$. 

At the point of catastrophe the assumptions of the implicit function theorem applied to the system of equations \eqref{tsar2} fail to hold true. Explicitly this means that
\beq\label{vanish1}
\det\left( a_{i,j}(\ub)\, t+A_{i,j}(\ub)\right)_{x=x_0, ~t=t_0, ~\ub=\ub_0}=\prod_{i=1}^n A_{i,i}(\ub_0)=0
\eeq
at the point of catastrophe.

\begin{defi} The point $x=x_0=0$, $t=t_0=0$, $\ub=\ub_0$ of gradient catastrophe is called {\rm generic} if at this point
\eqa\label{generic}
&&
A_{i,i}(\ub_0)=0, \quad A_{j,j}(\ub_0)\neq 0 \quad \mbox{\rm for}\quad j\neq i,
\nn\\
&&
\\
&&
A_{i,ii}(\ub_0)=0, \quad A_{i,iii}(\ub_0)\neq 0.
\nn
\eeqa
\end{defi}

The first part of this definition says that only $i$-th Riemann invariant breaks down at the point of catastrophe. The assumption $A_{i,ii}(\ub_0)=0$ is not a restriction; it says that the graph of the $i$-th Riemann invariant $u_i(x,0)$ has an inflection with vertical tangent at the point of catastrophe. As we will see the assumption $A_{i,iii}(\ub_0)\neq 0$ guarantees that
the inflection point does not degenerate.

Let us assume that $x=x_0=0$, $t=t_0=0$, $\ub=\ub_0=0$ is a generic point of gradient catastrophe with $i=n$. We will denote
$$
a_i^0:= a_i(\ub_0), \quad A_{i,i}^0:= A_{i,i}(\ub_0)
$$
etc. Observe that
$$
A_{i,j}^0=0 \quad \mbox{for all}\quad i\neq j
$$
due to the system \eqref{tsar1}.

Perform a linear change in the space of independent variables putting
\beq\label{change}
z_i:= x-a_i^0 t, \quad i=1, \dots, n.
\eeq

\begin{theorem} After the rescaling
\eqa\label{rescale1}
&&
z_m\mapsto k^{2/3} z_m, \quad m\neq n
\nn\\
&&
z_n\mapsto k\, z_n
\nn\\
&&
t\mapsto k^{2/3} t
\\
&&
u_m\mapsto k^{2/3} u_m, \quad m\neq n
\nn\\
&&
u_n\mapsto k^{1/3} u_n
\nn
\eeqa
the system \eqref{tsar1} converges to
\eqa\label{rescale2}
&&
z_m =A_{m,m}^0 u_m +{\mathcal O}\left( k^{1/3}\right), \quad m\neq n
\nn\\
&&
\\
&&
z_n =a_{n,n}^0 t\, u_n +\frac16 A_{n,nnn}^0 u_n^3 +{\mathcal O}\left(k^{1/3}\right).
\nn
\eeqa
\end{theorem}

{\sl Proof}. Spelling the equations \eqref{semiham2} out one obtains
\beq\label{tsar5}
a_{i,jk} =\frac{2 a_i -a_j -a_k}{(a_i-a_j)(a_i-a_k)} \, a_{i,j} a_{i,k}
+\frac{(a_i-a_j)^2 a_{i,k} a_{k,j} -(a_i-a_k)^2 a_{i,j} a_{j,k}}{(a_i-a_j)(a_i-a_k)(a_j-a_k)}
\eeq
for all triples $i$, $j$, $k$ of pairwise distinct indices. Using these equations along with the system \eqref{tsar1} yields

\begin{lemma} At the point of catastrophe one has
\eqa\label{ident0}
&&
A_{m,nn}^0=0, \quad A_{m,ni}^0=0\quad \mbox{\rm for}\quad m\neq n, \quad i\neq m, \, n 
\nn\\
&&
\nn\\
&&
A_{m,nnn}^0=0, \quad A_{m,nni}^0=0, \quad A_{m,nij}^0=0, \quad A_{m,ijk}^0=0 
\nn\\
&&
\mbox{\rm for}\quad m\neq n, \quad i, j, k \neq m, n, \quad i, j, k \quad \mbox{\rm distinct}
\nn\\
&&
\\
&&
A_{n,ni}^0=0, \quad A_{n,ij}^0=0, \quad \mbox{\rm for}\quad i, j\neq n, \quad i\neq j
\nn\\
&&
\nn\\
&&
A_{n,nni}^0=0, \quad A_{n,nij}^0=0, \quad A_{n,ijk}^0=0
\nn\\
&&
\mbox{\rm for}\quad i, j, k \neq n, \quad i, j, k \quad\mbox{\rm distinct}
\nn
\eeqa
\end{lemma}

\begin{cor} The Taylor expansions of the equations \eqref{tsar2}
near the point of catastrophe read
\eqa\label{taylor1}
&&
x-a_m^0t= t\,\left[ {\sum_i}' a_{m,i}^0 u_i +a_{m,m}^0 u_m+a_{m,n}^0 u_n +{\mathcal O}\left(|\ub|^2\right)\right] +A_{m,m}^0 u_m
\nn\\
&&
-\frac12\, {\sum_{i}}' \frac{a_{m,i}^0}{a_m^0 -a_i^0}\, A_{i,i}^0 u_i^2 + u_m\, {\sum_{i}}' \frac{a_{m,i}^0}{a_m^0 -a_i^0}\, A_{m,m}^0 u_i
\nn\\
&&
+\frac12\, A_{m,mm}^0 u_m^2 + u_m u_n \frac{a_{m,n}^0}{a_m^0 -a_n^0} \, A_{m,m}^0
\nn\\
&&
+\frac12\, {\sum_{i\neq j}}' A_{m,iij}^0 u_i^2 u_j +\frac16\, {\sum_i}' A_{m,iii}^0 u_i^3
\nn\\
&&
+\frac12\, u_m\, {\sum_{i\neq j}}' A_{m,mij}^0 u_i u_j +\frac12\, u_m {\sum_i}' A_{m,mii}^0 u_i^2
\nn\\
&&
+\frac12\, u_m^2{\sum_i}' A_{m,mmi}u_i +\frac16\, A_{m,mmm}^0 u_m^3 +{\mathcal O}\left(|\ub|^4\right).
\eeqa
\eqa\label{taylor2}
&&
x-a_n^0 t =t\, \left[ {\sum_i}' a_{n,i}^0 u_i +a_{n,n}^0 u_n +{\mathcal O}\left(|\ub|^2\right)\right]
-\frac12\, {\sum_i}' \frac{a_{n,i}^0}{a_n^0-a_i^0} \, A_{i,i}^0 u_i^2
\nn\\
&&
\\
&&
+\frac12\, {\sum_{i\neq j}}' A_{n,iij}^0 u_i^2 u_j +\frac16\, {\sum_i}' A_{n,iii}^0 u_i^3
+\frac12\, u_n {\sum_i}' A_{n,nii}^0 u_i^2 +\frac16\, A_{n,nnn}^0 u_n^3+{\mathcal O}\left(|\ub|^4\right).
\nn
\eeqa
\end{cor}

It is understood that the summation indices $i$, $j$ etc. in the sums denoted ${\sum}'$  are all distinct from $n$ and (in \eqref{taylor1}) from $m$.

Applying the rescaling \eqref{rescale1} to the equations \eqref{taylor1}, \eqref{taylor2} one completes the proof of the Theorem. \epf

\setcounter{equation}{0}
\setcounter{theorem}{0}
\section{Perturbations of nonlinear wave equations}\par

\subsection{Deformation problem-I: main definitions and the programme of research}\label{sect41}\par

Given a system of the first order quasilinear PDEs
\beq\label{def1}
\ub_t =A(\ub)\, \ub_x, \quad \ub=\left(u^1(x,t), \dots, u^n(x,t)\right)
\eeq
admitting a Hamiltonian description
\beq\label{def2}
\ub_t=\{ \ub(x), H_0\}_0, \quad H_0=\int h_0(\ub)\, dx
\eeq
with respect to a Poisson bracket of hydrodynamic type written in the flat coordinates in the form
\beq\label{def3}
\{ u^i(x), u^j(y)\}_0=\eta^{ij}\delta'(x-y), \quad \eta^{ji}=\eta^{ij}=\mbox{const}, \quad \det \left(\eta^{ij}\right) \neq 0,
\eeq
we say that a system of the form \eqref{pert0} is 
a {\it Hamiltonian deformation} of \eqref{def1} if it can be represented in the form
\beq\label{def4}
\ub_t =\{ \ub(x), H\}
\eeq
with a {\it perturbed Hamiltonian}
\eqa\label{def5}
&&
H=H_0 +\epsilon\, H_1 +\epsilon^2 H_2 +\dots
\nn\\
&&
\\
&&
H_k =\int h_k(\ub; \ub_x, \dots, \ub^{(k)})\, dx, \quad k\geq 1
\nn\\
&&
\nn\\
&&
\deg h_k(\ub; \ub_x, \dots, \ub^{(k)})=k
\nn
\eeqa
and a {\it perturbed Poisson bracket}
\eqa\label{def6}
&&
\{ u^i(x), u^j(y)\} =\{ u^i(x), u^j(y)\}_0+\epsilon\, \{ u^i(x), u^j(y)\}_1 +\epsilon^2 \{ u^i(x), u^j(y)\}_2+\dots
\nn\\
&&
\\
&&
\{ u^i(x), u^j(y)\}_k =\sum_{s=0}^{k+1} A^{ij}_{k\, s}\left(\ub(x);\ub_x(x), \dots, \ub^{(s)}(x) \right)\, \delta^{(k-s+1)}(x-y), \quad k\geq 1
\nn\\
&&
\nn\\
&&
\deg A^{ij}_{k\, s}\left(\ub;\ub_x, \dots, \ub^{(s)} \right)=s.
\nn
\eeqa
In the above formulae the notations like $f(\ub; \ub_x, \dots, \ub^{(k)})$ are used for differential polynomials
$$
f(\ub; \ub_x, \dots, \ub^{(k)})\in {\mathcal C}^\infty (M) [\ub_x, \dots, \ub^{(k)}], \quad \ub\in M,
$$
i.e., for polynomial functions on the jet bundle $J^k(M)$. The degree of a differential polynomial is defined by the rule
$$
\deg u^i_x=1, \quad \deg u^i_{xx}=2, \dots, \quad i=1, \dots, n.
$$
The sums in \eqref{def5} and \eqref{def6} are allowed to be infinite. However every term of the $\epsilon$-expansions of \eqref{def4} is a well-defined differential polynomial. Antisymmetry and Jacobi identity for the Poisson bracket are understood as identities for formal power series in $\epsilon$.
The delta-function symbol can be spelled out in the following way. Introduce a matrix of linear differential operators depending on the parameter $\epsilon$
\eqa\label{def7}
&&
\Pi^{ij} =\Pi^{ij}_0 +\epsilon\, \Pi^{ij}_1+\epsilon^2 \Pi^{ij}_2+\dots
\nn\\
&&
\\
&&
\Pi^{ij}_0 =\eta^{ij}\, \pal_x,
\nn\\
&&
\Pi^{ij}_k = \sum_{s=0}^{k+1} A^{ij}_{k\, s}\left(\ub;\ub_x, \dots, \ub^{(s)} \right)\, \pal_x^{k-s+1}, \quad k\geq 1.
\nn
\eeqa
One clearly has
\beq\label{def71}
\{ u^i(x), u^j(y)\} =\Pi^{ij}\, \delta(x-y).
\eeq
Then the perturbed Hamiltonian system reads
\beq\label{def8}
u^i_t =\Pi^{ij} \,\frac{\delta H}{\delta u^j(x)}
= \sum_{m\geq 0} \epsilon^m \sum_{k+l=m}\Pi^{ij}_k \, \frac{\delta H_l}{\delta u^j(x)}.
\eeq
The assumptions about the structure of the perturbative expansions in \eqref{def5} and \eqref{def6} imply that the terms of the expansion \eqref{pert0}
$$
B_m^i\left(\ub; \ub_x, \dots \ub^{(m+1)} \right) = \sum_{k+l=m}\Pi^{ij}_k \, \frac{\delta H_l}{\delta u^j(x)}, \quad m\geq 0, \quad i=1, \dots, n
$$
are differential polynomials of the degree $m+1$.

The class of Hamiltonian deformations is invariant with respect to (special) {\it Miura-type transformations} of the dependent variables
\beq\label{miur1}
\ub\mapsto \tilde \ub = \ub + \sum_{k\geq 1}\epsilon^k\, F_k\left(\ub; \ub_x, \dots, \ub^{(k)}\right) , \quad \deg F_k\left(\ub; \ub_x, \dots, \ub^{(k)}\right)=k.
\eeq
The transformation of the Hamiltonian is defined by the direct substitution; the Poisson bracket is transformed by the rule
\eqa\label{miur2}
&&
\{ \tilde u^i(x), \tilde u^j(y)\} =\tilde\Pi^{ij} \, \delta(x-y)
\nn\\
&&
\\
&&
\tilde\Pi^{ij} = L^i_p \,\Pi^{pq} {L^\dagger}^j_q
\nn
\eeqa
where $L$ and $L^\dagger$ are respectively the operator of linearization of the transformation \eqref{miur1} and its formal adjoint:
\beq\label{miur3}
{L}^i_k =\sum_s  \frac {\pal\tilde u^i}{\pal u^{k,s}} \pal_x^s, \quad
{L^\dagger}^i_k = \sum_s (-\pal_x)^s \, \frac {\pal\tilde u^i}{\pal u^{k,s}}.
\eeq

The transformations of the form \eqref{miur1} form a group. Two Hamiltonian deformations of the quasilinear system \eqref{def1} related by a transformation \eqref{miur1} are called {\it equivalent}. In particular the Hamiltonian deformation is called {\it trivial} if it is equivalent to the unperturbed system \eqref{def1}.

\begin{exam} The focusing/defocusing nonlinear Schr\"odinger equation
\beq\label{d-nls0}
i\,\epsilon\, \psi_t +\frac12\,\epsilon^2 \psi_{xx} \pm |\psi|^2\psi=0
\eeq
written in the coordinates
$$
u=|\psi|^2, \quad v=\frac{\epsilon}{2i}\left( \frac{\psi_x}{\psi} -\frac{\bar\psi_x}{\bar\psi}\right),
$$
takes the form of a Hamiltonian perturbation of the system \eqref{nlsd}:
\eqa\label{d-nls}
&&
u_t +(u\, v)_x=0
\\
&&
v_t +v\, v_x \mp u_x =\frac14 \epsilon^2\left( \frac{u_{xx}}{u}-\frac12\, \frac{u_x^2}{u^2}\right)_x.
\nn
\eeqa
The Poisson bracket gets no perturbation:
$$
\{ u(x), v(y)\} =\delta'(x-y), \quad \{ u(x), u(y)\}=\{v(x), v(y)\}=0.
$$
The perturbed Hamiltonian reads
\beq\label{d-nls1}
H=H_0 +\epsilon^2\, H_2, \quad H_0 =\int \frac12\, (u^2 \mp u\, v^2)\, dx, \quad H_2= -\int \frac{u_x^2}{8\, u} \, \, dx.
\eeq
\end{exam}

\begin{exam} The generalized FPU system with the Hamiltonian \eqref{fpu1} in the translation invariant variables
\eqa
&&
w_n = q_n -q_{n-1}
\nn\\
&&
v_n=p_n
\eeqa
and interpolation
\eqa
&&
w_n=w(\epsilon \, n, \epsilon\, t)
\nn\\
&&
v_n=v(\epsilon \, n, \epsilon\, t)
\nn
\eeqa
yields the following equations for the functions $w=w(x,t)$, $v=v(x,t)$
\eqa\label{d-fpu1}
&&
w_t(x) = \frac1{\epsilon}\,\left[v(x)-v(x-\epsilon)\right]=v_x -\frac12\epsilon\, v_{xx} +\frac16\epsilon^2 v_{xxx} +\dots
\nn\\
&&
\\
&&
v_t(x) =\frac1{\epsilon}\,\left[P'(w(x+\epsilon)) -P'(w(x))\right]
=\pal_xP'(w) +\frac12 \epsilon\, \pal_x^2 P'(w)+\frac16 \epsilon^2 \,\pal_x^3 P'(w)+\dots.
\nn
\eeqa
This system is a Hamiltonian perturbation of the nonlinear wave equation \eqref{wave2}. The Poisson brackets of the variables $w(x)$, $v(x)$ have a nontrivial dependence on $\epsilon$. Indeed, the brackets of the discrete variables read
$$
 \{ w_n, v_m\} =\delta_{n,m} -\delta_{n, m+1},
$$
other brackets vanish.
After interpolation and division by $\epsilon$ one obtains the perturbed bracket
\eqa\label{d-fpu2}
&&
\{ w(x), v(y)\} =\frac1{\epsilon}\,\left[ \delta(x-y) - \delta(x-y-\epsilon)\right]
\nn\\
&&
\\
&&
\quad\qquad\qquad\,\,\,=\delta'(x-y)-\frac12\,\epsilon\,\delta''(x-y) +\frac16\, \epsilon^2 \delta'''(x-y) +\dots.
\nn
\eeqa
One can reduce the Poisson bracket \eqref{d-fpu2} to the standard form
$$
\{ u(x), v(y)\} =\delta'(x-y)
$$
by the Miura-type transformation \eqref{uw}. Indeed, introducing the shift operator
$$
\Lambda=e^{\epsilon\, \pal_x}
$$
rewrite the bracket \eqref{d-fpu2} as follows
$$
\{ w(x), v(y)\} =\frac1{\epsilon} \left( 1-\Lambda^{-1}\right) \,\delta(x-y).
$$
The substitution \eqref{uw} represented in the form 
$$
u=\frac{\epsilon\, \pal_x}{1-\Lambda^{-1}}\, w
$$
implies
$$
\{ u(x), v(y)\} = \frac{\epsilon\, \pal_x}{1-\Lambda^{-1}}\, \{ w(x), v(y)\} = \pal_x\delta(x-y).
$$
In the new variables the deformed system
is described by the perturbed Hamiltonian:
\eqa
&&
u_t =\pal_x \frac{\delta H}{\delta v(x)}
\nn\\
&&
\nn\\
&&
v_t =\pal_x \frac{\delta H}{\delta u(x)}
\nn
\eeqa
\eqa\label{hami4}
&&
H=\int h\, dx=\int \left[\frac12\, v^2(x) +P\left( w(x) -w(x-\epsilon)\right)\right]\, dx
\nn\\
&&
\\
&&
h=\frac12\, v^2 + P(u) -\frac{\epsilon^2}{24} \, P''(u)\, u_x^2
+\frac{\epsilon^4}{5760}\,\left[ 8\, P''(u)\, u_{xx}^2 - P^{IV}(u) \, u_x^4\right] + {\mathcal O}(\epsilon^6)
\nn
\eeqa
(modulo inessential total derivatives). The series in $\epsilon$ in this case are all infinite.
\end{exam}

In fact the following general statement about deformations of the Poisson bracket  \eqref{def3} holds true.

\begin{theorem} 1). Any deformation \eqref{def5} of the Poisson bracket \eqref{def3} is trivial.

\noindent 2). Any special Miura-type transformation \eqref{miur1} canonical  with respect to the bracket $\{ ~,~\}=\{~, ~\}_0$
$$
\{ \tilde u^i(x), \tilde u^j(y)\} = \{ u^i(x), u^j(y)\}
$$
is the time-$\epsilon$ shift generated by a Hamiltonian $\Phi$:
\eqa\label{time-shift}
&&
\tilde u^i= u^i + \epsilon\, \{u^i(x), \Phi\} +\frac{\epsilon^2}2\, \{\{u^i(x), \Phi\} , \Phi\} +\dots
\nn\\
&&
\\
&&
\Phi = \int\varphi \, dx
\nn\\
&&
\varphi= \sum_{k\geq 0}\epsilon^k \varphi_k(\ub; \ub_x, \dots, \ub^{(k)}), \quad \deg \varphi_k(\ub; \ub_x, \dots, \ub^{(k)})=k, \quad k\geq 0
\nn
\eeqa
with some differential polynomials $\varphi_0(\ub)$, $\varphi_1(\ub; \ub_x)$, \dots.
\end{theorem}

\pf follows from triviality in positive degrees of Poisson cohomology of the bracket \eqref{def3} established in \cite{getzler}, \cite{magri} (see details in \cite{DZ}). 

Thanks to the Theorem one reduces the study of Hamiltonian deformations of the Hamiltonian quasilinear equations to the study of systems written in the standard form
\beq\label{stand-ham1}
u^i_t =\eta^{ij} \pal_x \frac{\delta H}{\delta u^j(x)}, \quad i=1, \dots, n
\eeq
just perturbing the Hamiltonian,
\eqa\label{hpert1}
&&
H=H^{[0]} + \epsilon \, H^{[1]} +\epsilon^2 H^{[2]} +\dots
\nn\\
&&
\\
&&
H^{[k]} =\int h^{[k]}(\ub; \ub_x, \dots, \ub^{(k)}) \, dx, \quad 
\deg h^{[k]}(\ub; \ub_x, \dots, \ub^{(k)})=k, \quad k\geq 0
\nn
\eeqa
with some differential polynomials $h^{[0]}(\ub)$, $h^{[1]}(\ub; \ub_x)$, $h^{[2]}(\ub; \ub_x, \ub_{xx})$, \dots. Here 
$$
H^{[0]} =\int h^{[0]}(\ub)\, dx
$$ 
is the Hamiltonian of the unperturbed equation. The equivalence of perturbed Hamiltonians is established by the canonical transformations \eqref{time-shift}:
\beq\label{time-shift1}
H\mapsto \tilde H = H +\epsilon\, \{ H,\Phi\} +\frac{\epsilon^2}2\, \{\{ H,\Phi\}, \Phi\} +\dots.
\eeq

In the subsequent sections we will discuss the main problems of the deformation theory of Hamiltonian PDEs, namely,

{\bf Problem 1}: describe first integrals/infinitesimal symmetries of a given deformed Hamiltonian system by studying the procedure of {\it extension} of integrals/symmetries of the unperturbed system to integrals/symmetries of the perturbed system.

{\bf Problem 2}: classify {\it integrable} Hamiltonian perturbations for which {\it all} first integrals/infinitesimal symmetries of the unperturbed system can be extended to first integrals/infinitesimal symmetries of the perturbed one.

{\bf Problem 3}: comparative study of the {\it properties of solutions} of the perturbed versus unperturbed systems. In a more specific way, characterize the perturbations for which the solutions to the perturbed Cauchy problem, in a suitable class of initial data, exists for the times beyond the moment of gradient catastrophe for the solution to the unperturbed Cauchy problem. For such a class of perturbed Hamiltonian systems
establish the {\it universal properties} of solutions near the point
of gradient catastrophe of the unperturbed system, that is, those properties independent from the choice of the initial data. A more ambitious project is in establishing universality independently on the choice of Hamiltonian perturbation.

\subsection{Deformation problem-II: extension of infinitesimal symmetries and $D$-operator}\par

Let
\eqa\label{do1}
&&
H_{\rm pert}=H^{[0]} +\epsilon\, H^{[1]}+\epsilon^2 H^{[2]}+\dots
\nn\\
&&
H^{[k]}=\int h^{[k]}(\ub; \ub_x, \dots, \ub^{(k)})\, dx, \quad \deg h^{[k]}(\ub; \ub_x, \dots, \ub^{(k)})=k, \quad k\geq 0
\\
&&
h^{[0]}=\frac12\, v^2 +P(u)
\nn
\eeqa
be a perturbation of the Hamiltonian of the nonlinear wave equation. We know that, given a solution $f=f(u,v)$ to the linear PDE \eqref{chap1} the Hamiltonians $H^{[0]}$ and
\beq\label{do2}
H_f^{[0]}=\int f\, dx
\eeq
commute:
$$
\{ H^{[0]}, H_f^{[0]}\}=0.
$$
The goal is to construct a deformation
\eqa\label{dof}
&&
H_f=H_f^{[0]} +\epsilon\, H_f^{[1]}+\epsilon^2 H_f^{[2]}+\dots
\nn\\
&&
H_f^{[k]}=\int h_f^{[k]}(\ub; \ub_x, \dots, \ub^{(k)})\, dx, \quad \deg h_f^{[k]}(\ub; \ub_x, \dots, \ub^{(k)})=k, \quad k\geq 0
\\
&&
h_f^{[0]}=f(u,v)
\nn
\eeqa
for an arbitrary solution $f$ of \eqref{chap1} such that 
$$
\{ H_{\rm pert}, H_f\}=0.
$$

\begin{defi} We say that the linear differential operator
\eqa\label{doper}
&&
D=D^{[0]} +\epsilon\, D^{[1]} +\epsilon^2 D^{[2]}+\dots
\nn\\
&&
\nn\\
&&
D^{[0]}={\rm id}, \quad D^{[k]} =\sum b^{[k]}_{i_1 \dots i_{m(k)}} (\ub; \ub_x, \dots, \ub^{(k)} )\frac{\pal^{m(k)}}{\pal u^{i_1} \dots \pal u^{i_{m(k)}}}
\\
&&
\nn\\
&&
\deg b^{[k]}_{i_1 \dots i_{m(k)}} (\ub; \ub_x, \dots, \ub^{(k)} )=k, \quad k\geq 1
\nn
\eeqa
is a $D$-{\rm operator} for the perturbation \eqref{do1} if for any solution $f$ to \eqref{chap1} the Hamiltonian
\beq\label{hdf}
H_f =\int D f\, dx
\eeq
commutes with $H_{\rm pert}$.
\end{defi}

In the formula \eqref{doper} $m(k)$ is some positive integer depending on $k$. It is easy to see that
\beq\label{3k}
m(k)=\left[ \frac{3k}2\right].
\eeq
The summation is taken over all indices $i_1$, \dots, $i_{m(k)}$ from 1 to 2. As usual the coefficients $b^{[k]}_{i_1 \dots i_{m(k)}} (\ub; \ub_x, \dots, \ub^{(k)} )$ are differential polynomials. As the $D$-operator makes sense only acting on the kernel of
\beq\label{2oper}
\pal_u^2 -P''(u) \pal_v^2,
\eeq
the coefficients are not determined uniquely. One can choose a suitable normal form for this operator, e.g., keeping only derivatives of the form
$$
\pal_v^m \quad \mbox{and}\quad \pal_u\pal_v^{m-1}.
$$

We will show that the $D$-operator is determined uniquely just from the commutativity with $H_{\rm pert}$, as it follows from the following

\begin{theorem} \label{hyp} The  $D$-operator, if exists, is determined uniquely up to a total $x$-derivative, i.e., if $D_1$ and $D_2$ are two $D$-operators then
$$
D_2f - D_1f\in \mbox{\rm Im}\, \pal_x \quad \forall\, f\in \mbox{\rm Ker} (\pal_u^2 -P''(u)\, \pal_v^2).
$$
\end{theorem}

\pf  follows from the following result saying that the nonlinear wave equation has no first integrals depending nontrivially on the derivatives.

\begin{lemma} \label{lm-hyp} Let $P(u)$ be a potential such that $P''(u)\neq 0$. Let
$$
g=g(\ub; \ub_x, \dots, \ub^{(k)})
$$
be a differential polynomial of the degree $k>0$ such that the functional
$$
G=\int g(\ub; \ub_x, \dots, \ub^{(k)})\, dx
$$
commutes with  $H^{[0]}$,
\beq\label{doper2}
\{ G, H^{[0]}\} =0 
\eeq
Then $g$ is a total $x$-derivative:
$$
g(\ub; \ub_x, \dots, \ub^{(k)})=\pal_x\tilde g(\ub; \ub_x, \dots, \ub^{(k-1)}).
$$
\end{lemma}

The proof of this lemma is straightforward and will appear elsewhere.

\begin{cor} Let the perturbed Hamiltonian \eqref{do1} admit a $D$-operator.
Then for any pair of solutions $f$, $g$ to \eqref{chap1} the Hamiltonians $H_f$ and $H_g$ commute:
\eqa\label{doper1}
&&
\{ H_f, H_g\} =0, \quad H_f =\int D f\, dx, \quad H_g=\int D g\, dx
\nn\\
&&
\\
&&
f_{uu}=P''(u) f_{vv}, \quad g_{uu}=P''(u) g_{vv}.
\nn
\eeqa
\end{cor}

\pf From Jacobi identity it follows that 
$$
\{\{ H_f, H_g\}, H_{\rm pert}\}=0.
$$
Assume the Poisson bracket $\{ H_f, H_g\}$ to have the form
$$
\{ H_f, H_g\}=\epsilon^k F^{[k]} +{\mathcal O}\left( \epsilon^{k+1}\right)
$$
for some $k\geq 1$. The functional $F^{[k]}$ must commute with $H^{[0]}$. Due to the Lemma it must vanish. \epf.

\begin{remark} One can consider a more general class of the first order quasilinear Hamiltonian systems written in the form
\eqa\label{gener1}
&&
u_t =\pal_x h_v
\nn\\
&&
v_t=\pal_x h_u
\eeqa
for a given smooth function $h=h(u,v)$. The conservation laws $f=f(u,v)$ for this system are determined from the linear PDE
\beq\label{gener2}
h_{vv} f_{uu} = h_{uu} f_{vv}.
\eeq
The associated Hamiltonians
$$
H_f^{[0]}=\int f\, dx
$$
commute pairwise. We say that the system \eqref{gener1} is {\rm regular}\footnote{Similarly to the definition of regular elements in Lie algebras.} if it has no nontrivial local conservation laws depending on higher derivatives (cf. Lemma \ref{lm-hyp}). It would be interesting to characterize all regular quasilinear Hamiltonian systems.
\end{remark}

We end this Section with considering a simple example of Hamiltonian perturbations of the {\it scalar} transport equation
$$
u_t+u\, u_x=0
$$
classified in \cite{du2}. The Hamiltonian structure of the transport equation reads
\eqa
&&
u_t+\{ u(x), H^{[0]}\} \equiv u_t +\pal_x \frac{\delta H^{[0]}}{\delta u(x)}=0
\nn\\
&&
H^{[0]}=\frac16 \int u^3 \, dx, \quad \{ u(x), u(y)\}=\delta'(x-y).
\nn
\eeqa
All perturbations of the Hamiltonian $H^{[0]}$ up to the order $\epsilon^4$ are parameterized by two arbitrary functions of one variable $c=c(u)$ and $p=p(u)$ as follows:
\beq
H_{\rm pert}=\int \left[ \frac{u^3}6 - \frac{\epsilon^2}{24} c(u)\, u_x^2
+\epsilon^4  p(u) \,u_{xx}^2 \right]\, dx.
\nn
\eeq
This perturbation turns out to be {\it integrable} with the order $\epsilon^4$ approximation
$$
\{ H_f, H_g\} = {\mathcal O}\left( \epsilon^6\right), \quad H_f=\int D_{c,p}f\, dx, \quad H_g =\int D_{c,p}g\, dx, \quad \forall f=f(u), ~g=g(u).
$$
The associated $D$-operator has the form
\eqa\label{old}
&&
D_{c,p}f = f-\frac{\epsilon^2}{24} c\, f''' u_x^2+\epsilon^4\left[\left( p\,f''' + 
     \frac{{c}^2\,f^{(4)}}{480} \right) \,u_{xx}^2\right.
\nn\\
&&
\\
&&\left.  - 
 \left(   
      \frac{c\,c''\,f^{(4)}}{1152} + \frac{c\,c'\,f^{(5)}}{1152} + 
     \frac{{c}^2\,f^{(6)}}{3456} +\frac{p'\,f^{(4)}}{6} +
     \frac{p\,f^{(5)}}{6}\right)\, u_{x}^4\right] +{\mathcal O}\left( \epsilon^6\right),
\nn
\eeqa
the function $f$ in this case is arbitrary.

\subsection{Deformation problem-III: an integrability test}\par

Existence of a $D$-operator for a given perturbed Hamiltonian gives rise to a simple integrability test, as we will now explain for the example of the generalized FPU system. The perturbed Hamiltonian reads
\beq\label{itest1}
H_{\rm pert}=\int\left[\frac12\, v^2 + P(u) -\frac{\epsilon^2}{24} \, P''(u)\, u_x^2 \right]\, dx +{\mathcal O}\left(\epsilon^4\right)
\eeq
(see \eqref{hami4} above). Under what condition for the potential $P(u)$ a $D$-operator exists within the $\epsilon^2$ approximation? We will show that this is possible only under the constraint
\beq\label{itest2}
P'' P^{IV} =(P''')^2.
\eeq
To derive \eqref{itest2} let us look for the perturbed Hamiltonian density with indetermined coefficients
\eqa
&&
H_f =\int h_f\, dx, \quad \{ H_f, H_{\rm pert}\}={\mathcal O}\left(\epsilon^3\right)
\nn\\
&&
h_f =f + \epsilon\, h_f^{[1]}(\ub; \ub_x) +\epsilon^2 h_f^{[2]}(\ub; \ub_x, \ub_{xx}).
\nn
\eeqa
The first correction must be linear in $u_x$, $v_x$. Adding if necessary a total $x$-derivative one can reduce to the study of the first perturbation of the form
$$
h_f^{[1]}=p(u,v)\, v_x
$$
for some function $p=p(u,v)$.
Computation of the bracket gives
$$
\{ H_f, H_{\rm pert}\} = \epsilon\, \int p_u\,\left[ v_x^2 -P''(u) u_x^2\right]\, dx+{\mathcal O}\left(\epsilon^2\right).
$$
It is easy to see that the integrand is never a total $x$-derivative unless $p_u=0$, i.e. $p=p(v)$. Hence $h_f^{[1]}$ must be itself a total $x$-derivative.

Let us now consider the order 2 terms. Up to a total $x$-derivative they can be written in the form
$$
h_f^{[2]}= \frac12 \left(a(u,v) u_x^2 + 2 b(u,v) u_x v_x + c(u,v) v_x^2\right).
$$
Computation of the order $\epsilon^2$ terms in the Poisson bracket gives the following functional:
\eqa
&&
\{H_f, H_{\rm pert}\}=
\nn\\
&&
=\epsilon^2 \int \left\{\frac1{12} f_v P''u_{xxx}+
\left[ \left( \frac16 f_v P''' - b\, P''\right) \, u_x -a\, v_x\right]\, u_{xx} -\left( c\, P'' u_x +b\, v_x\right)\, v_{xx}\right.
\nn\\
&&
+\frac1{24}\left( f_v P^{IV} +12 P'' a_v -24 P'' b_u\right) \, u_x^3
-\frac12\left( a_u +2 P'' c_u\right)\, u_x^2 v_x -\frac12\left( 2a_v +P'' c_v\right)\, u_x v_x^2
\nn\\
&&
\left.
+\frac12\left( c_u-2 b_v\right)\, v_x^3\right\}\, dx+{\mathcal O}\left(\epsilon^4\right).
\nn
\eeqa
Denote $I$ the integrand. One must verify the equations
$$
{\rm E}_u I=0, \quad {\rm E}_v I=0.
$$
Vanishing of the coefficient of $u_{xxx}$ in ${\rm E}_v I$ gives
$$
a=\left(c-\frac1{12} f_{vv} \right)P''.
$$
Next, vanishing of the coefficient of $v_{xx} v_x$ in ${\rm E}_v I$ implies
$$
b_v=c_u.
$$
Hence it exists a potential $\lambda=\lambda(u,v)$ such that
$$
b=\lambda_u, \quad c=\lambda_v.
$$
After the substitution we obtain, collecting the coefficients of $u_{xx} u_x$, $u_{xx} v_x$, $v_{xx} u_x$, $v_{xx} v_x$, the following system of equations
\eqa
&&
6 P''' \lambda_v + P'' f_{uvv}=0
\nn\\
&&
P''' \lambda_u -P'' \lambda_{uu} +{P''}^2 \lambda_{vv}=0
\nn\\
&&
6 \lambda_{uu} -6P''\lambda_{vv} +P'' f_{vvv}=0.
\nn
\eeqa
It follows that
\eqa
&&
\lambda_u =-\frac16\frac{{P''}^2}{P'''} f_{vvv}
\nn\\
&&
\lambda_v=-\frac16 \frac{P''}{P'''} f_{uvv}.
\nn
\eeqa
Equating the mixed derivatives
$$
\left( \lambda_u\right)_v =\left( \lambda_v\right)_u
$$
implies
$$
\frac{(P''')^2 -P'' P^{IV}}{6\, (P''')^2} f_{uvv}=0.
$$
This proves \eqref{itest2}.

\medskip

The equation \eqref{itest2} for the potential can be easily solved: dividing by $P''' P''$ yields
$$
\frac{P'''}{P''} = \frac{P^{IV}}{P'''}.
$$
Integration gives
$$
P'''=c\, P''.
$$
Hence
$$
P(u)=k\, e^{c\, u} +a\, u + b
$$
for some constants $a$, $b$, $c$, $k$.

We have proved that the only generalized FPU system passing the integrability test must have exponential potential. This gives Toda lattice.

\begin{remark} It would be interesting to compare the above arguments with the theory of small amplitude normal forms
for finite FPU chains developed recently by A.Henrici and T.Kappeler \cite{hk1, hk2}.
\end{remark}

\subsection{Deformation problem-IV: on classification of integrable perturbations}\par

Let $f=f(u,v)$ and $g=g(u,v)$ be two first integrals of the nonlinear wave equation \eqref{wave2}, i.e., two solutions to the linear PDE \eqref{chap1}. The Hamiltonians
$$
H_f^{[0]}=\int f\, dx, \quad H_g^{[0]}=\int g\, dx
$$
commute:
$$
\{ H_f^{[0]}, H_g^{[0]}\}=0.
$$
We obtain a commutative Lie algebra of Hamiltonians parameterized by solutions to \eqref{chap1}. The goal is to classify deformations of this algebra of the form \eqref{def4}
preserving the commutativity.

So far the classification is available only within the $\epsilon^2$
approximation. Let us describe this preliminary result and outline the classification procedure.

\begin{theorem} Any 2-integrable perturbation of the nonlinear wave equation \eqref{wave2} is parameterized by two arbitrary functions of one variable $\rho_\pm =\rho_\pm (r_\pm)$. The associated $D$-operator reads
\eqa\label{psi-d}
&&
D_{\rho, r} f= f +\frac{\epsilon^2}2 \left\{ \left[P''( \rho_u f_{vvv} + \rho_v f_{uvv}) +\frac12 P''' \rho_v f_{vv}\right]\, u_x^2
\right.
\nn\\
&&
+2\left( P''\rho_v f_{vvv} +\rho_u f_{uvv} +\frac{P'''}{4 P''} \rho_u f_{vv}\right)\, u_x v_x
\nn\\
&&
\left.
+\left( \rho_u f_{vvv} +\rho_v f_{uvv}\right) \, v_x^2
\right\}+{\mathcal O}( \epsilon^3),
\nn
\eeqa
where
$$
\rho=\rho(u,v)=\rho_+(r_+) +\rho_-(r_-),
$$
the Riemann invariants $r_\pm=r_\pm(u,v)$ are defined in \eqref{riem2}.
\end{theorem}

\pf Let $f=f(u,v)$, $g=g(u,v)$ be two arbitrary solutions to the linear PDE \eqref{chap1}. Any ${\mathcal O}(\epsilon)$-perturbation of the commuting Hamiltonians $H_f^{[0]}$ and $H_g^{[0]}$ can be written in the form
$$
H_f =H_f^{[0]} +\epsilon\int p(u,v) v_x \, dx+{\mathcal O}\left(\epsilon^2\right), \quad H_g =H_g^{[0]} +\epsilon\int q(u,v) v_x \, dx+{\mathcal O}\left(\epsilon^2\right)
$$
for some functions $p$ and $q$. From the commutativity
$$
\{ H_f, H_g\}={\mathcal O}\left(\epsilon^2\right)
$$
it follows, calculating the coefficient of $u_{xx}$ in
$\frac{\delta}{\delta u(x)} \{ H_f, H_g\}$, that
\beq
\frac{\pal_u p}{f_{vv}} = \frac{\pal_u q}{g_{vv}}.
\eeq
Consider now a canonical transformation
$$
\ub \mapsto \ub +\epsilon \{ \ub(x), K\}+{\mathcal O}\left(\epsilon^2\right)
$$
generated by the Hamiltonian
$$
K=\int k(u,v)\, dx.
$$
The Hamiltonians $H_f$ and $H_g$ will transform according to the following rule:
$$
H_f \mapsto H_f +\epsilon \{ H_f^{[0]}, K\} +{\mathcal O}\left(\epsilon^2\right), \quad H_g \mapsto H_g +\epsilon \{ H_g^{[0]}, K\} +{\mathcal O}\left(\epsilon^2\right).
$$
Choosing the function $k=k(u,v)$ satisfying
$$
k_{uu}-P''(u) k_{vv} = \frac{\pal_u p}{f_{vv}} = \frac{\pal_u q}{g_{vv}}
$$
one eliminates the terms linear in $\epsilon$ from both perturbed Hamiltonians $H_f$ and $H_g$.

At the next step one has to consider the perturbation of order $\epsilon^2$ that can be chosen in the following form:
\eqa
&&
H_f =H_f^{[0]} +\frac{\epsilon^2}2 \int \left( a_f u_x^2 + 2 b_f u_x v_x + c_f v_x^2\right) \, dx +{\mathcal O}\left(\epsilon^3\right)
\nn\\ 
&&
H_g =H_g^{[0]} +\frac{\epsilon^2}2 \int \left( a_g u_x^2 + 2 b_g u_x v_x + c_g v_x^2\right) \, dx +{\mathcal O}\left(\epsilon^3\right)
\nn
\eeqa
with some functions $a_i=a_i(u,v)$, $b_i=b_i(u,v)$, $c_i=c_i(u,v)$, $i=f,\, g$.  From the commutativity
$$
\{ H_f, H_g\}={\mathcal O}\left(\epsilon^3\right)
$$
it follows, calculating the coefficient of $v_{xxx}$ in $\frac{\delta}{\delta u(x)} \{ H_f, H_g\}$, that
$$
\left( a_f -c_f P''(u)\right) g_{vv} = \left( a_g -c_g P''(u)\right) f_{vv}.
$$
Hence
\eqa
&&
a_f = c_f P''(u) + k(u,v) f_{vv}
\nn\\
&&
a_g =c_g P''(u) + k(u,v) g_{vv}
\nn
\eeqa
for some function $k(u,v)$.

At the next step we substitute
\eqa
&&
b_f = r_1 f_{vvv} + r_2 f_{uvv} +r_3 f_{vv} + r_4 f_{uv}, \quad c_f = s_1 f_{vvv} + s_2 f_{uvv} +s_3 f_{vv} + s_4 f_{uv}
\nn\\
&&
b_g = r_1 g_{vvv} + r_2 g_{uvv} +r_3 g_{vv} + r_4 g_{uv}, \quad c_g = s_1 g_{vvv} + s_2 g_{uvv} +s_3 g_{vv} + s_4 g_{uv}
\nn
\eeqa
with some coefficients $r_i=r_i(u,v)$, $s_i=s_i(u,v)$, $i=1, \dots, 4$ (the coefficients of the $D$-operator). The calculations similar to above yield the following expressions for the coefficients:
\eqa
&&
r_1 =P'' \rho_v, \quad r_2 =\rho_u, \quad r_3 =\frac{P'''}{4 P''}\, \rho_u, \quad r_4=0
\nn\\
&&
s_1=\rho_u, \quad s_2 =\rho_v, \quad s_4=0
\nn\\
&&
k=\frac12 P''' \rho_v -2 s_3\, P''
\nn
\eeqa
where the function $s_3=s_3(u,v)$ is arbitrary and the function $\rho=\rho(u,v)$ satisfies the linear PDE
$$
\rho_{uu} -P'' \rho_{vv} =\frac{P'''}{2P''} \,\rho_u.
$$
It is readily seen that this equation is equivalent to the D'Alembert equation in Riemann invariants:
$$
\rho_{uu} -P'' \rho_{vv} -\frac{P'''}{2P''} \,\rho_u =-4 P'' \frac{\pal^2 \rho}{\pal r_+ \pal r_-}.
$$
Finaly, we obtain the
$\epsilon^2$-correction of the Hamiltonian densities in the following form
\eqa
&&
h_f = f +\frac{\epsilon^2}2 \left\{ \left[P''( \rho_u f_{vvv} + \rho_v f_{uvv}) +\frac12 P''' \rho_v f_{vv}\right]\, u_x^2
\right.
\nn\\
&&
+2\left( P''\rho_v f_{vvv} +\rho_u f_{uvv} +\frac{P'''}{4 P''} \rho_u f_{vv}\right)\, u_x v_x
\nn\\
&&
\left.
+\left( \rho_u f_{vvv} +\rho_v f_{uvv}\right) \, v_x^2
+ s_3 \,  \left( v_x^2 -P'' u_x^2\right)\, f_{vv}\right\}+{\mathcal O}\left( \epsilon^3\right),
\nn
\eeqa
a similar expression holds for $h_g$ replacing $f\mapsto g$. It is easy to see that the term containing the arbitrary function $s_3$ can be eliminated by a canonical transformation generated by the Hamiltonian
$$
\Phi=\frac{\epsilon}2 \int \phi(u,v) v_x\, dx
$$
with $\phi=\phi(u,v)$ determined from
$$
\phi_u =s_3.
$$ \epf

\subsection{Quasitriviality and solutions to the perturbed
nonlinear wave equations}\par

Recall (see Section \ref{sect41} above) that the $\epsilon$-perturbation of a first order quasilinear system is called trivial if it can be eliminated by a (special) Miura-type transformation \eqref{miur1}. All terms in the $\epsilon$-expansion of this transformation are differential polynomials. 

\begin{defi} We say that the perturbation is {\rm quasitrivial} if it can be eliminated by a transformation of the form
\beq\label{qmiur1}
\ub\mapsto \tilde \ub = \ub + \sum_{k\geq 1}\epsilon^k\, F_k\left(\ub; \ub_x, \dots, \ub^{(m(k))}\right) , \quad \deg F_k\left(\ub; \ub_x, \dots, \ub^{(m(k))}\right)=k
\eeq
where the terms $F_k$ for $k\geq 1$ are {\rm rational} functions in the jet variables $\ub_x$, $\ub_{xx}$, \dots $\ub^{(m(k))}$ for some integer $m(k)$.
\end{defi}

As it was proved in \cite{DLZ}, all {\it bihamiltonian} perturbations are quasitrivial. Moreover, the order $m(k)$ of the highest
derivative in $F_k$ can be estimated as
\beq\label{3g-2}
m(k) =\left[\frac{3k}2\right].
\eeq
It would be interesting to extend the Quasitriviality Theorem of \cite{DLZ} to all Hamiltonian perturbations\footnote{The quasitriviality theorem was proved by S.-Q.Liu and Y.Zhang in \cite{LZ2} for a more general class of perturbations of scalar PDEs. In the recent paper \cite{LWZ} the quasitriviality technique was applied to the analysis of existence of a local Hamiltonian structure for scalar evolutionary PDEs.}. 

Existence of a quasitriviality transformation gives a universal formula for the perturbative solutions to the deformed systems
with very general classes of initial data. Namely, substituting a solution $\ub_0=\ub_0 (x,t)$ of the unperturbed system into \eqref{qmiur1} one obtains a formal series in $\epsilon$ satisfying the deformed PDE in all orders. In the examples of Section \ref{sect5} below all the rational functions in the quasitriviality transformation have the denominators that are powers of
\beq\label{det-den}
\Delta=v_x^2-P''(u) u_x^2.
\eeq

In the Hamiltonian case the quasitriviality transformation can be represented as a time-$\epsilon$ canonical transformation (see \cite{du2}). That means that there exists a Hamiltonian
\eqa\label{qmiur2}
&&
K =\sum_{m\geq 1} \epsilon^m K_m
\nn\\
&&
\\
&&
K_m =\int k_m(\ub; \ub_x, \dots, \ub^{()})\, dx, \quad m\geq 1
\nn
\eeqa
such that the canonical transformation
\beq\label{qmiur3}
\ub\mapsto \tilde\ub=\ub+ \epsilon\{ \ub(x), K\} +\frac{\epsilon^2}2 \{\{ \ub(x), K\} , K\} +\dots
\eeq
transforms {\it all} unperturbed commuting Hamiltonians into the perturbed ones,
\beq\label{qmiur4}
H_f =\int f(\tilde \ub)\, dx \quad \mbox{for any}\quad f\quad\mbox{satisfying}\quad f_{uu}=P''(u) f_{vv}.
\eeq
The generating Hamiltonian $K$ has to be found from the equation
\beq\label{qmiur5}
H_f^{[0]} +\epsilon\{ H_f^{[0]}, K\} +\frac{\epsilon^2}2\{ \{ H_f^{[0]}, K\}, K\} +\dots =H_f \equiv H_f^{[0]}+\epsilon\, H_f^{[1]}+\epsilon^2 H_f^{[2]}+\dots.
\eeq
The Hamiltonians constructed in this way will automatically commute; the highly nontrivial set of cancellations
ensures that the densities of these Hamiltonians up to total $x$-derivatives are {\it polynomials} in the jet coordinates $\ub_x$, $\ub_{xx}$ etc.

The above perturbative expansions can be used for asymptotic integration of the deformed Hamiltonian systems for sufficiently small times (although rigorous analytic results justifying validity of these asymptotic expansions are still to be established). They clearly lose their validity when approaching the time of gradient catastrophe since the denominators become all equal to zero. In order to study the behaviour of the perturbed solutions near the critical point
we will now describe an alternative procedure of constructing solutions to the perturbed Hamiltonian systems. Let us first prove the following

\begin{lemma} Let the perturbed Hamiltonian  of the nonlinear wave equation have the form
\beq\label{v-ind}
H_{\rm pert}=\int \left[\frac12\, v^2 + P(u) +\sum_{k\geq 1}\epsilon^k h^{[k]}(u; u_x, v_x, \dots, u^{(k)}, v^{(k)})\right]\, dx
\eeq
(i.e., no explicit dependence on $v$ in the higher corrections). The the $D$-operator, if exists, has $v$-independent coefficients.
\end{lemma}

\pf Indeed, under the above assumptions the coefficients of the commutativity equation $\{ H_f, H_{\rm pert}\}=0$ for $H_f =\int Df\, dx$
\eqa
&&
{\rm E}_u \left[ \frac{\delta H_f}{\delta u(x)} \left( v_x + \epsilon\pal_x\frac{\delta \tilde H_{\rm pert}}{\delta v(x)}\right) 
+ \frac{\delta H_f}{\delta v(x)} \left( P''(u)u_x +\epsilon\pal_x \frac{\delta \tilde H_{\rm pert}}{\delta u(x)}\right) \right]=0
\nn\\
&&
{\rm E}_v \left[ \frac{\delta H_f}{\delta u(x)} \left( v_x + \epsilon\pal_x\frac{\delta \tilde H_{\rm pert}}{\delta v(x)}\right) 
+ \frac{\delta H_f}{\delta v(x)} \left( P''(u) u_x+\epsilon\pal_x \frac{\delta \tilde H_{\rm pert}}{\delta u(x)}\right) \right]=0
\nn
\eeqa
do not depend on $v$ (see \eqref{kommut}). Here
$$
\tilde H_{\rm pert} = \int\sum_{k\geq 1} \epsilon^{k-1} h^{[k]}(u; u_x, v_x, \dots, u^{(k)}, v^{(k)})\, dx.
$$
Moreover, the shift 
$$
v\mapsto v+c
$$
maps to itself the space of solutions to the linear PDE \eqref{chap1}. Therefore the same shift applied to the $D$-operator produces another $D$-operator. Due to the uniqueness theorem \ref{hyp} the shifted operator coincides with $D$ after adding, if necessary, a total $x$-derivative. \epf 

\begin{remark} The class of Hamiltonan perturbations considered in the Lemma corresponds to the following subclass
of perturbed nonlinear wave equations:
\beq\label{v-ind1}
u_{tt} =P''(u)u_{xx} +\epsilon\pal_x a(u; u_x, u_t, u_{xx}, u_{xt}, \dots; \epsilon).
\eeq
\end{remark}

Given a solution $h=h(u,v)$ to the PDE \eqref{chap1}, consider a Hamiltonian perturbation of the quasilinear system
\eqa
&&
u_t =\pal_x h_v
\nn\\
&&
v_t =\pal_x h_u
\nn
\eeqa
(we have changed notation $s\mapsto t$ for the time variable).
Solutions to this system are parameterized in the form \eqref{h-eq} by functions $f=f(u,v)$ satisfying the same PDE \eqref{chap1} according to the method of characteristics. Let us assume that the perturbed Hamiltonian equation admits a $v$-independent $D$-operator. We will show that the solutions to the perturbed equation close to \eqref{h-eq} can be obtained from a system of  ``infinite order ODEs" written in a variational form. This system of ODEs will be called {\it string equation}.

\begin{theorem} The solutions to the 
{\rm string equation}
\eqa\label{string}
&&
x+ t \,\frac{\delta H_{h'}}{\delta u(x)} =\frac{\delta H_{f}}{\delta u(x)}
\nn\\
&&
\\
&&
 \qquad t \,\frac{\delta H_{h'}}{\delta v(x)} =\frac{\delta H_{f}}{\delta v(x)}
\nn
\eeqa
also solve the Hamiltonian equations
\eqa\label{fl1}
&&
 u_t =\pal_x \frac{\delta H_{h}}{\delta v(x)} 
 \nn\\
 &&
 \\
 &&
 v_t =\pal_x \frac{\delta H_{h}}{\delta u(x)} 
\nn
\eeqa
 where $f=f(u,v)$ is another solution to $f_{uu}=P''(u) f_{vv}$,
 $$
 h':= \frac{\pal h}{\pal v}.
 $$
\end{theorem}

\pf Since the coefficients of $D$-operator are $v$-independent, one has an identity
\beq\label{v'1}
 \frac{\pal}{\pal v} Dh = Dh'.
\eeq
Due to obvious commutativity
$$
\left[ \frac{\pal}{\pal v}, {\rm E}_u\right] =\left[ \frac{\pal}{\pal v}, {\rm E}_v\right]=0
$$
in the assumptions of the Theorem from \eqref{v'1} one derives
\beq\label{v'2}
 \frac{\pal}{\pal v} \frac{\delta H_h}{\delta u(x)} =\frac{\delta H_{h'}}{\delta u(x)}, \quad \frac{\pal}{\pal v} \frac{\delta H_h}{\delta v(x)} =\frac{\delta H_{h'}}{\delta v(x)}.
\eeq
Let us use the above identities in order to prove that the flow
\eqa\label{gal1}
&&
u_s = \qquad t\, \pal_x \frac{\delta H_{h'}}{\delta v(x)}
\nn\\
&&
\\
&&
v_s =1+ t\, \pal_x \frac{\delta H_{h'}}{\delta u(x)}
\nn
\eeqa
is an infinitesimal symmetry of the system \eqref{fl1}. Indeed,
\eqa
&&
(u_t)_s= t\, \{ \pal_x \frac{\delta H_h}{\delta v(x)}, H_{h'}\} +\pal_x \frac{\pal}{\pal v} \frac{\delta H_h}{\delta v(x)} = t\, \{ \{ u(x), H_h\}, H_{h'}\} + \pal_x \frac{\delta H_{h'}}{\delta v(x)}
\nn\\
&&
\nn\\
&&
(v_t)_s = t\, \{ \{ v(x), H_h\}, H_{h'}\} + \pal_x \frac{\delta H_{h'}}{\delta u(x)}.
\nn
\eeqa
Computation of the derivatives in the opposite order yields
\eqa
&&
(u_s)_t =\pal_x \frac{\delta H_{h'}}{\delta v(x)} + t\, \{ \pal_x \frac{\delta H_{h'}}{\delta v(x)}, H_h\} = \pal_x \frac{\delta H_{h'}}{\delta v(x)} +\{ \{ u(x), H_{h'}\}, H_h\}
\nn\\
&&
\nn\\
&&
(v_s)_t = \pal_x \frac{\delta H_{h'}}{\delta u(x)} +\{ \{ v(x), H_{h'}\}, H_h\}.
\nn
\eeqa
The Jacobi identity and commutativity of the Hamiltonians
$\{ H_h, H_{h'}\}=0$ imply
\eqa
&&
\{ \{ u(x), H_h\}, H_{h'}\} =\{ \{ u(x), H_{h'}\}, H_h\}
\nn\\
&&
\{ \{ v(x), H_h\}, H_{h'}\} =\{ \{ v(x), H_{h'}\}, H_h\}.
\nn
\eeqa
So
$$
(u_t)_s=(u_s)_t, \quad (v_t)_s=(v_s)_t.
$$

The commutativity remain valid adding to the symmetry \eqref{gal1} an arbitrary flow of the hierarchy, e.g., replacing \eqref{gal1} with
\eqa\label{gal2}
&&
u_s = \qquad t\, \pal_x \frac{\delta H_{h'}}{\delta v(x)}-\pal_x \frac{\delta H_{f}}{\delta v(x)}
\nn\\
&&
\\
&&
v_s =1+ t\, \pal_x \frac{\delta H_{h'}}{\delta u(x)}-\pal_x \frac{\delta H_{f}}{\delta u(x)}
\nn
\eeqa
for some solution $f$  to the linear PDE \eqref{chap1}. It is wellknown (see, e.g., \cite{n74, lax}) that the set of stationary points of an infinitesmal symmetry is invariant for the flow \eqref{fl1}. Integrating one obtains the string equation \eqref{string}. In the dispersion limit the string equation coincides with the equation \eqref{hodo2} of the method of characteristics:
\eqa
&&
x+t\, h_{uv} = f_{u}
\nn\\
&&
\qquad t\, h_{vv}=f_{v}.
\nn
\eeqa

\subsection{Special Painlev\'e-type transcendents and Universality Conjectures}\label{sect46}\par

We have seen in Section \ref{sect3} above that locally near the point of catastrophe all generic solutions to a nonlinear wave equation and, more generally, to any Hamiltonian PDE commuting with the wave equation, have a very standard behaviour. Namely, if the catastrophe takes place in the domain of hyperbolicity then the local behaviour is described by the Whitney singularity equations. For the catastrophe in the domain of ellipticity the local behaviour is given in terms of elliptic umbilic singularity. In both cases the shape of the solution is locally universal, i.e., independent from the choice of the solution, up to a simple affine transformation of the independent variables and a nonlinear transformation of the dependent variables.

The Universality Conjecture for the perturbed systems says that, similarly, the leading term in the asymptotic expansion of the solution at $\epsilon\to 0$ near a critical point does not depend on the choice of a generic solution, up to simple transformations of the independent/dependent variables. The universal shapes that occur in such a local description are given in terms of certain very particular solutions to Painlev\'e equations and their generalizations. Let us first describe these solutions. The description will be given separately for catastrophes taking place in the domain of hyperbolicity or domain of ellipticity.

$\bullet$ Type I singularity. Consider the following 4th order ODE 
for the function $U=U(X)$ depending on the parameter $T$
\beq\label{P12}
X=T\, U -\left[ \frac16 U^3 +\frac1{24} ( {U'}^2 + 2 U\, U'' ) +\frac1{240} U^{IV}\right], \quad U'=\frac{dU}{dX} \quad \mbox{ etc.}
\eeq
The equation \eqref{ode0} can be considered as a higher order analogue of the Painlev\'e-I equation. It is often identified as the $P_I^2$ equation. 
As it was proved by T.Claeys and M.Vanlessen \cite{cl1}, this equation has a particular solution without singularities on the real axis $X\in\mathbb R$ for all real values of the parameter $T$. Denote this solution $U(X,T)$. We refer to the paper \cite{cl1} for the precise characterization of the particular solution in terms of the Stokes multipliers of the associated with \eqref{ode0} linear differential operator with polynomial coefficients.
For large $|X|$ the solution grows like
$$
U(X,T)\sim (-6\, X)^{1/3}.
$$
It is monotone decreasing for $T<< 0$; for positive $T$ it develops a zone of oscillations (see Fig.1).

\vspace{.2in}
\centerline {
\includegraphics[width=2.6in]{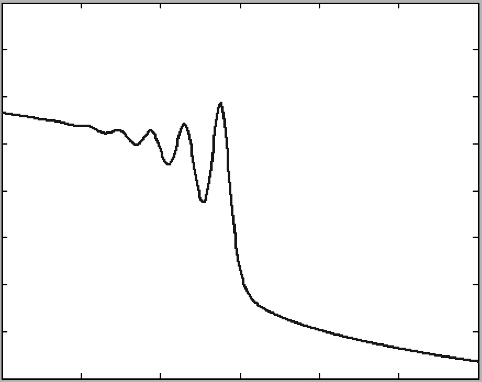}
}
\vspace{.2in}
\begin{center} Fig.1: The solution to \eqref{ode0} for $T>0$\end{center}

\smallskip

We are ready to formulate the Universality Conjecture describing the local behaviour of solutions to the perturbed Hamiltonian system \eqref{fl1} near the Type I critical point
$(x_{\rm c}, t_{\rm c}, u_{\rm c}, v_{\rm c})$. Denote $r_\pm$ the Riemann invariants \eqref{riem2} of the unperturbed system. Let
$$
x_\pm =(x-x_{\rm c}) + c_\pm (t-t_{\rm c}), \quad c_\pm=h_u(u_{\rm c}, v_{\rm c}) \pm \sqrt{P''(u_{\rm c})} h_v(u_{\rm c}, v_{\rm c})
$$
be the characteristic variables at the point of catastrophe. Let us assume that at the point of catastrophe of the unperturbed system the {\it second}  Riemann invariant $r_-$ breaks down, while the invariant $r_-$ remains smooth (see \eqref{tip1} above).
Then
one should expect the following asymptotic behaviour of solutions to \eqref{fl1} near the critical point:
\eqa\label{univer1}
&&
r_+ = r_+ (u_{\rm c}, v_{\rm c}) + c_+ x_+ +a_+\, \epsilon^{4/7} U''\left( b_- \frac{x_-}{\epsilon^{6/7}}, b_+\frac{x_+}{\epsilon^{4/7}}\right)
+{\mathcal O}\left(\epsilon^{6/7}\right)
\nn\\
&&
\\
&&
r_-=r_- (u_{\rm c}, v_{\rm c})+a_- \epsilon^{2/7} U\left( b_- \frac{x_-}{\epsilon^{6/7}}, b_+\frac{x_+}{\epsilon^{4/7}}\right) +{\mathcal O}\left(\epsilon^{4/7}\right)
\nn
\eeqa
for some constants $a_\pm$, $b_\pm$, $c_+$ depending on the choice of solution and on the perturbed Hamiltonian. Here $U=U(X,T)$ is the particular solution to the 4th order ODE \eqref{ode0} described above.

In the next Section we will explain how to compute the constants for the particular examples of integrable systems.

$\bullet$ Type II singularity. Here we need to specify a particular solution $W=W_0(Z)$ to the classical Painlev\'e-I ($P_I$) equation
\beq\label{p1}
W_{ZZ} =6\, W^2 -Z.
\eeq
The particular solution in question was discovered by Boutroux
in the beginning of the last century \cite{bo}. It is characterized by the asymptotic behaviour
$$
W\sim -\sqrt{\frac{Z}6}, \quad |\arg Z|<\frac{4\pi}5, \quad |Z|\to\infty
$$
(the principal branch of the square root is chosen).
In particular the solution $W_0(Z)$ is free of poles within the sector $|\arg Z|<\frac{4\pi}5$ for sufficiently large $|Z|$. Boutroux called it the {\it tritronqu\'ee} solution. A.Kapaev \cite{ka} gave a complete characterization of the {\it tritronqu\'ee} solution in terms of particular values of the Stokes multipliers
of the linear differential operator with polynomial coefficients associated with \eqref{p1}.

\vspace{.2in}
\centerline {
\includegraphics[width=5in]{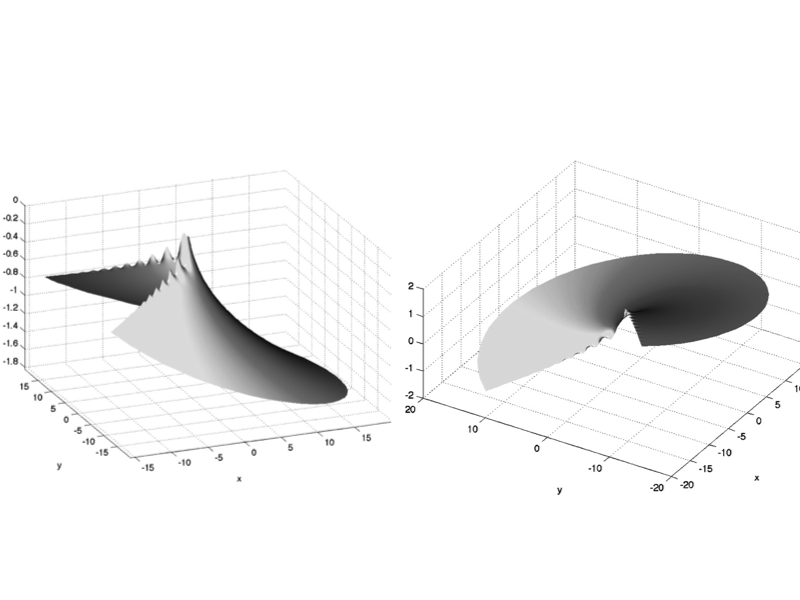}
}
\vspace{.2in}
\begin{center} Fig.2: The real (left) and imaginary parts of the {\it tritronqu\'ee} solution to Painlev\'e-I equation in the sector $|\arg Z|<\frac{4\pi}5$.\end{center}

In \cite{dgk} it was conjectured that {\it all} poles of the {\it tritronqu\'ee} solution to $P_I$ are confined within the sector
$|\arg Z|\geq \frac{4\pi}5$. The Fig. 2 borrowed from \cite{dgk}
shows the graphs of the real and imaginary parts of the {\it tritronqu\'ee} solution suggesting that, indeed there are no poles in the sector $|\arg Z|< \frac{4\pi}5$.

The Universality Conjecture for the Type II critical behaviour suggests the following structure of the leading term of the asymptotic behaviour of (complex) Riemann invariants
\beq\label{univer2}
r_+ = r_+(u_{\rm c}, v_{\rm c}) + \beta\, \epsilon^{2/5} \, W_0\left( \frac{x_+}{\alpha\, \epsilon^{4/5}}\right) +{\mathcal O}\left( \epsilon^{4/5}\right),
\eeq
the dependence of the complex conjugate Riemann invariant $r_-$ on the variable $x_-$ is given by complex conjugation.
Here $W_0(Z)$ is the {\it tritronqu\'ee} solution to $P_I$,
$$
x_+ = (x-x_{\rm c}) +(h_v^0 +i\, \sqrt{-P''_0}\, h_v^0)\, (t-t_{\rm c}),
$$
$\alpha$, $\beta$ are some complex numbers. These numbers are such that for any suffuiciently small $|t|$ and $\epsilon >0$ the line
$$
Z= \frac{x_+}{\alpha\, \epsilon^{4/5}}, \quad x\in \mathbb R
$$
on the complex plane belongs entirely to the sector
$$
|\arg Z| <\frac{4\pi}5
$$
where the {\it tritronqu\'ee} solution $W_0(Z)$ conjecturally has no poles.

In the next Section we will consider in more details the most widely known examples of Boussinesq, Toda, NLS and Ablowitz - Ladik equations where are approach can be applied. For the case of Boussinesq equation we will give more detailed explanations that motivate our Universality Conjectures. We will see that, although these equations look so different, the local structure of critical behaviour in all these cases is given by the same formulae. The rigorous mathematical justification of such Universality is still to have been done; a more challenging problem is to establish Universality for a more general class of Hamiltonian perturbations, without restricting to only integrable perturbations.

\smallskip

\setcounter{equation}{0}
\setcounter{theorem}{0}
\section{Main examples}\label{sect5}\par

\subsection{Boussinesq equation}\par

As it was already explained in Section \ref{sect21}, the Boussinesq equation can be considered as a Hamiltonian perturbation of the nonlinear wave equation with the potential
$$
P(u)=-\frac16 u^3.
$$ 
The linear PDE \eqref{chap1} for this potential 
coincides with the Tricomi equation
\beq\label{tric}
f_{uu}+u\, f_{vv}=0.
\eeq
The domain of hyperbolicity is the half-plane $u<0$.
The $D$-operator transforming solutions of this equation to first integrals for Boussinesq hierarchy reads
\eqa\label{d-bous}
&&
f\mapsto D_{\rm Bous}f = f +\epsilon^2\left[ \left(\frac12\, f_{vv} + u\, f_{uvv}\right) u_x^2 +2 u\, f_{vvv} u_x v_x -f_{uvv} v_x^2\right]
\\
&&
+\epsilon^4 \left[ \left( \frac65\, u^2 f_{4v} -f_{uvv}\right) u_{xx}^2 -\frac45\left( 3u\, f_{u\, 3v} + 2 f_{vvv}\right) u_{xx} v_{xx}
-\frac65\, u\, f_{4v} v_{xx}^2\right.
\nn\\
&&
-\frac15\, f_{4v} u_{xx} v_x^2 -f_{u\, 3v} v_{xx} u_x^2-\frac1{120}\left( 31\, f_{4v} -20 u^3 f_{6v} + 92 u\, f_{u\, 4v}\right) \, u_x^4
\nn\\
&&
-\frac2{15}\, u \left( 13\, f_{5v} +5u\, f_{u\, 5v}\right)\, u_x^3 v_x+\frac1{10}\left( 11\, f_{u\, 4v} -10u^2 f_{6v}\right)\, u_x^2 v_x^2
\nn\\
&&\left.
+\frac1{15} \left( 7\, f_{5v} + 10u\, f_{u\, 5v}\right)\, u_x v_x^3+\frac16 \, u\, f_{6v} v_x^4
\right]+{\mathcal O}\left(\epsilon^6\right).
\nn
\eeqa
Here and below we adopt the following notations for high order partial derivatives
$$
f_{4v} := \frac{\pal^4 f}{\pal v^4}, \quad f_{u\, 3v}:=\frac{\pal^4 f}{\pal u \,\pal v^3}
$$
etc. In particular, the choice
$$
f=\frac12\, v^2 -\frac16\, u^3
$$
gives the Hamiltonian density of Boussinesq equation
$$
D_{\rm Bous} f = \frac12\, v^2 -\frac16\, u^3 +\frac{\epsilon^2}2 u_x^2.
$$

Let us proceed to the study of critical behaviour of solutions to an equation of the Boussinesq hierarchy
\eqa\label{bous_hi1}
&&
u_t =\pal_x \frac{\delta H}{\delta v(x)}
\nn\\
&&
\\
&&
v_t =\pal_x \frac{\delta H}{\delta u(x)}
\nn
\eeqa
where
\beq\label{bous_hi2}
H=\int D_{\rm Bous} \tilde h\, dx, \quad \tilde h_{uu}+u \,\tilde h_{vv}=0.
\eeq

We begin with Type I critical points; for them $u_{\rm c} <0$. It is convenient to change the sign
$$
u\mapsto -u.
$$
Consider the string equation
\beq\label{uras}
\delta \int \left[x\, u + t\,  D_{\rm Bous} h -D_{\rm Bous} f\right]\, dx=0
\eeq
for arbitrary two solutions $f=f(u,v)$, $h=h(u,v)$ to 
$$
f_{uu}=u\, f_{vv}, \quad h_{uu} = u\, h_{vv}.
$$
Here the solution $h=h(u,v)$ is related with the Hamiltonian
\eqref{bous_hi2} of the system \eqref{bous_hi1} (with $u$ replaced by $-u$) by
\beq\label{bous_hi3}
\pal_v \tilde h = h.
\eeq
At the critical point $(u_{\rm c}, v_{\rm c})$ one must have
\eqa
&&
f_u^0=f_v^0=0
\nn\\
&&
f_{uv}^0 =\sqrt{u_{\rm c}} f_{vv}^0
\nn\\
&&
f_{uvv}^0 =-\frac{f_{vv}^0}{4\, u_{\rm c}} + \sqrt{u_{\rm c}}\, f_{vvv}^0
\nn
\eeqa
(see \eqref{crit11} and \eqref{generic14} above). Let us expand the string equation near the critical point in the Taylor series with respect to the Riemann invariants coordinates
$$
r_\pm = v\pm \frac23 u^{3/2}.
$$
It is also convenient to perform a linear change of independent variables $(x,t)\mapsto (x_+, x_-)$ 
$$
x_\pm = x+\left( h_u^0 \pm \sqrt{u_{\rm c}} \,h_v^0\right)\, t
$$
(cf. \eqref{xpm}). After the rescaling
\eqa\label{scal-hyp}
&&
\bar x_+ \mapsto k^{2/3} \bar x_+, \quad \bar x_- \mapsto k\, \bar x_-
\nn\\
&&
\\
&&
\bar r_+ \mapsto k^{2/3} \bar r_+, \quad \bar r_-\mapsto k^{1/3} \bar r_-
\nn\\
&&
\nn\\
&&
\epsilon\mapsto k^{7/6} \epsilon
\nn
\eeqa
of the shifted variables
$$
\bar x_\pm = x_\pm -x_\pm ^0, \quad \bar r_\pm =r_\pm -r_\pm^0= \bar v\pm \sqrt{u_{\rm c}}\, \bar u
$$
at the limit $k\to 0$ one arrives at the following system
\eqa\label{strun1}
&&
\bar x_+ = 2\sqrt{u_{\rm c}} \, f_{vv}^0\,\bar r_+ +\frac{\epsilon^2}{2\sqrt{u_{\rm c}}}f_{vv}^0 \,\pal_-^2 \bar r_- +{\mathcal O}\left( k^{1/3}\right)
\\
&&
\nn\\
&&\label{strun2}
\bar x_-=A_0 \bar x_+ \bar r_- +B_0 \bar r_-^3 +\epsilon^2 A_2 \left[ \left( \pal_- \bar r_-\right)^2 +2\bar r_-\, \pal_-^2 \bar r_-\right]
\\
&&
-\frac{\epsilon^2}{2 \sqrt{u_{\rm c}}} f_{vv}^0\, \pal_-^2 \bar r_+
+\epsilon^4 \tilde A_4\, \pal_-^4 \bar r_- +{\mathcal O}\left( k^{1/3}\right)
\nn
\eeqa
with some coefficients $A_0$, $B_0$, $A_2$, $\tilde A_4$ (see below). In the above equations we denote
$$
\pal_\pm :=\frac{\pal}{\pal x_\pm}.
$$
From the first equation one has
$$
\pal_-^2 \bar r_+ = -\frac{\epsilon^2}{4 u_{\rm c}}\, \pal_-^4 \bar r_-+{\mathcal O}\left( k^{1/3}\right).
$$
Substitution in \eqref{strun2} yields
\beq\label{strun3}
\bar x_-=A_0 \bar x_+ \bar r_- +B_0 \bar r_-^3 +\epsilon^2 A_2 \left[ \left( \pal_- \bar r_-\right)^2 +2\bar r_-\, \pal_-^2 \bar r_-\right]
+\epsilon^4  A_4\, \pal_-^4 \bar r_- +{\mathcal O}\left( k^{1/3}\right)
\eeq
where
\eqa
&&
A_0 =\frac{h_{vv}^0}{2 h_v^0} -\frac{h_{uv}^0}{2 \sqrt{u_{\rm c}}\, h_v^0} -\frac1{8u_{\rm c}^{3/2}}
\nn\\
&&
\nn\\
&&
B_0 =\frac16\, \left[ \frac{5 f_{vv}^0}{32 u_{\rm c}^{5/2}} +\frac1{4u_{\rm c}} f_{vvv}^0 -\sqrt{u_{\rm c}}\, f_{vvvv}^0 + f_{uvvv}^0\right]
\nn\\
&&
\nn\\
&&
A_2 = 6\, \sqrt{u_{\rm c}}\, B_0
\nn\\
&&
\nn\\
&&
A_4=\frac{72}{5} u_{\rm c}\, B_0.
\nn
\eeqa
The substitution
\eqa
&&
\bar x_- = \alpha\, X
\nn\\
&&
\bar x_+ =\gamma\, T
\nn\\
&&
\bar r_- =\beta\, U
\nn
\eeqa
with the choice
$$
k=\epsilon^{6/7}
$$
and
\eqa
&&
\alpha=-2^{10/7} 3^{4/7} u_{\rm c}^{3/14} B_0^{1/7}
\nn\\
&&
\nn\\
&&
\beta= \left(\frac23\right)^{1/7} u_{\rm c}^{1/14} B_0^{-2/7}
\nn\\
&&
\nn\\
&&
\gamma=-2^{9/7} 3^{5/7} u_{\rm c}^{1/7} A_0^{-1}B_0^{3/7} 
\nn
\eeqa
reduces the equation \eqref{strun3} to the standard form
$$
X= U\, T -\left[ \frac16 U^3 + \frac1{24}\left({U'}^2 + 2 U\, U''\right) +\frac1{240} U^{IV}\right] +{\mathcal O}\left( \epsilon^{2/7}\right).
$$
The equation \eqref{strun1} takes the following form
$$
\bar x_+ =2 \sqrt{u_{\rm c}} f_{vv}^0 \left( \bar r_+ + \frac1{4 u_{\rm c}} \frac{\beta}{\alpha}^2 U''\right) +{\mathcal O}\left( \epsilon^{2/7}\right).
$$
Returning to the original variables yields the following final form of the Universality Conjecture for the Type I critical behaviour for generic solutions to any equation to the  Boussinesq hierarchy
\eqa\label{konec1}
&&
r_+ = r_+^0 + \frac1{2 \sqrt{u_{\rm c}} f_{vv}^0 }
\left[ x-x_{\rm c} + \left( h_u^0 +\sqrt{u_{\rm c}} h_v^0\right)\, (t-t_{\rm c})\right] 
\nn\\
&&
\nn\\
&&
-\frac{\epsilon^{4/7} }{4 u_{\rm c}} \frac{\beta}{\alpha^2} U''\left( \frac{x-x_{\rm c} + \left( h_u^0 -\sqrt{u_{\rm c}} h_v^0\right) (t-t_{\rm c})}{\epsilon^{6/7}\alpha}, \frac{x-x_{\rm c} + \left( h_u^0 +\sqrt{u_{\rm c}} h_v^0\right) (t-t_{\rm c})}{\epsilon^{4/7}\gamma}\right) +{\mathcal O}\left( \epsilon^{6/7}\right)
\nn\\
&&
\\
&&
r_- =r_-^0 + \beta\, \epsilon^{2/7} U\left( \frac{x-x_{\rm c} + \left( h_u^0- \sqrt{u_{\rm c}} h_v^0\right) (t-t_{\rm c})}{\epsilon^{6/7}\alpha}, \frac{x-x_{\rm c} + \left( h_u^0 +\sqrt{u_{\rm c}} h_v^0\right) (t-t_{\rm c})}{\epsilon^{4/7}\gamma}\right) +{\mathcal O}\left( \epsilon^{4/7}\right).
\nn
\eeqa
Here the function $U(X,T)$ smooth for all real values $X$, $T$ is determined from the $P_I^2$ equation, as it was explained above in Section \ref{sect46}. For the particular choice
$$
h(u,v)=v
$$
one obtains the critical behaviour of solutions to the Boussinesq equation itself. The solutions are parameterized by functions $f=f(u,v)$ satisfying
$$
f_{uu}=u\, f_{vv}.
$$

Let us now proceed to considering the Type II critical behaviour.
In this case it is convenient to use the original choice of signs for Boussinesq equation. At the critical point $u_{\rm c}>0$. Given a pair of solutions $f=f(u,v)$, $h=h(u,v)$ to Tricomi equation \eqref{tric}, consider the string equation of the form \eqref{uras}. At the Type II critical point one has
$$
f_u^0=f_v^0=0, \quad f_{uu}^0 =f_{uv}^0=f_{vv}^0=0.
$$
The Riemann invariants
$$
r_\pm =v\mp \frac23\, i\, u^{3/2}
$$
are complex conjugate. For the unperturbed system they are holomorphic resp. antiholomorphic functions of the characteristic variables
$$
x_\pm =x+(h_u \mp i\, \sqrt{u}\, h_v)t,
$$
i.e.
$$
\frac{\pal r_+}{\pal x_-} = \frac{\pal r_-}{\pal x_+}=0.
$$
The Universality Conjecture for the Type II critical behaviour first formulated for the case of focusing NLS by T.Grava, C.Klein and the author \cite{dgk} says that the holomorphicity remains valid also for the leading term of asymptotic expansion of the perturbed solution. Moreover it was conjectured in \cite{dgk} that this leading term is given in terms of the {\it tritronqu\'ee} solution to the Painlev\'e-I equation (see above Section \ref{sect46}). Let us reproduce the arguments of \cite{dgk} adapted for the case of Type II critical behaviour of solutions to Boussinesq equation (i.e., choosing $h(u,v)=v$).

As above we expand the string equations in Taylor series near the critical point. After the rescaling 
\eqa
&&
\bar x_\pm \mapsto k\, \bar x_\pm
\nn\\
&&
\bar r_\pm \mapsto k^{1/2} \bar r_\pm
\nn\\
&&
\epsilon\mapsto k^{5/4}\epsilon
\nn
\eeqa
of the shifted variables 
\eqa
&&
\bar x_\pm = (x-x_{\rm c}) \pm i\, \sqrt{u_{\rm c}}\, (t-t_{\rm c})
\nn\\
&&
\bar r_\pm =r_\pm - r_\pm (u_{\rm c}, v_{\rm c})=\bar v \pm i\, \sqrt{u_{\rm c}} \bar u
\nn
\eeqa
one obtains at $k\to 0$
\eqa\label{prep1}
&&
\bar x_+ = \lambda \, (\bar  r_+^2 -\epsilon^2 \mu\, \bar r_+'')+{\mathcal O}\left( k^{1/2}\right)
\\
&&
\lambda= \frac12 ( f_{uvv}^0 -i\, \sqrt{u_{\rm c}}\, f_{vvv}^0), \quad \mu=4i \sqrt{u_{\rm c}}.
\nn
\eeqa
A complex conjugate equation holds valid for the variables $x_-$, $r_-$. 

In the above equation
$$
\bar r_+'' =\pal_x^2 \bar r_+ =\pal_+^2 \bar r_+ +2\pal_+ \pal_- \bar r_+ + \pal_-^2 \bar r_+.
$$
Performing the same rescaling in the Boussinesq equation we find
\eqa
&&
\bar u_t =\bar v_x
\nn\\
&&
\bar v_t + u_{\rm c} \bar u_x ={\mathcal O}\left( k^{1/2}\right).
\nn
\eeqa
This is equivalent to
$$
\pal_- \bar r_+ ={\mathcal O}\left(k^{1/2}\right).
$$
So, at the leading order we may assume that
$\bar r_+$ is an analytic function in $x_+$. The equation \eqref{prep1} then becomes
\beq\label{prep2}
\bar x_+ = \lambda \, (\bar  r_+^2 -\epsilon^2 \mu\, \pal_+^2\bar r_+)+{\mathcal O}\left( k^{1/2}\right).
\eeq
Choosing
$$
k=\epsilon^{4/5}
$$
and doing a change of variables
\eqa
&&
x_+ = \alpha\, Z
\nn\\
&&
\bar r_+ = \beta \, W
\nn
\eeqa
with
\eqa
&&
\alpha= (6 \lambda\, \mu^2)^{1/5}
\nn\\
&&
\beta= \left( 6^3 \frac{\mu}{\lambda^2}\right)^{1/5}
\nn
\eeqa
we arrive at the needed equation
$$
Z=6\, W^2 - W'' +{\mathcal O}\left( \epsilon^{2/5}\right)
$$
(see Remark 5.2 in \cite{dgk} explaining the choice of the fifth root of the complex numbers). Returning to the original variables one obtains the Universality Conjecture for the Type II critical behaviour for equations of Boussinesq hierarchy in the form
\eqa\label{bous-ell}
&&
v-v_{\rm c}-i\, u_{\rm c}^{1/2}\left( u- u_{\rm c}\right) \nn\\
&&
\\
&&
=\epsilon^{2/5} \beta\, W_0\left( \frac{x-x_{\rm c} +(h_u^0 -i \sqrt{u_{\rm c}}\, h_v^0)(t-t_{\rm c})}{\alpha\, \epsilon^{4/5}} \right) 
+{\mathcal O}\left( \epsilon^{4/5}\right)
\nn
\eeqa
where $W_0(Z)$ must be the {\it tritronqu\'ee} solution to the Painlev\'e-I equation. Indeed, it is easy to see that the line
$$
Z=\frac1{\alpha} \left( x-x_{\rm c} +(h_u^0 -i \sqrt{u_{\rm c}}\, h_v^0)(t-t_{\rm c})\right), \quad x\in \mathbb R
$$
belongs to the sector
$$
|\arg Z| <\frac{4\pi}5.
$$

\subsection{Toda hierarchy and its extension}\par

The linear PDE \eqref{chap1} is of hyperbolic type for all real $(u,v)$:
\beq\label{toda-pde}
f_{uu} = e^u\, f_{vv}.
\eeq
The $D$-operator transforming solutions of this equation into conserved quantities for (extended; see \cite{cdz}) Toda hierarchy reads
\eqa
&&
f\mapsto h_f=D_{\rm Toda} f
\nn\\
&&
=f-\frac{ {\epsilon}^2}{24}\,\left[ 
         \left(f_{vv} + 
         2\,f_{uvv} \right)\,v_x^2 +2\,e^u\,\left(f_{uvv}\,
          u_x^2 + 
         2\,f_{vvv}\,u_x\,v_x\right)  \right]   + 
  {\epsilon}^4\,\left\{ 
     \frac{1}{720}\,f_{vv}\,v_{xx}^2 
         \right.\nn\\
         &&\left.
          + 
     \frac{1}{120}e^u\,f_{vvv}\,
        \left( u_x^2 + u_{xx} \right) \,
        v_{xx}+ 
     \frac{1 }{4320} f_{uvv}\,\left[ 2\,e^u\,\left(u_x^4 - 
          3u_{xx}^2\right) + 
          3\,u_{xx}\,v_x^2 + 24\,v_{xx}^2 \
\right]
\right.\nn\\
          && \left.
            + 
     \frac{1}{4320}f_{uvvv}\,
        \left[ 4\,e^u\,\left(2u_x^3\,v_x  + 
          15\,u_x^2\,v_{xx} + 
          18\,u_{xx}\,v_{xx}\right) - 
          3\,u_x\,v_x^3\right] 
\right.
\nn\\
&&\left.
 + \frac{1}{5760} f_{vvvv}\,
        \left[ 16\,e^{2\,u}\,\left( 
          3\,u_{xx}^2 -u_x^4 \right)- 4\,e^u\,\left(
          7\,u_x^2\,v_x^2 +8\,u_{xx}\,v_x^2  -12\,v_{xx}^2\right)- 
          v_x^4 \right] 
                    \right.\nn\\
          &&\left.
- \frac{1}{2160}e^u\,f_{5v}\,u_x\,
        v_x\,\left( 14\,e^u\,u_x^2 + 
          15\,v_x^2 \right)- 
     \frac{1}{4320} f_{u\,4v}\,
        \left( 17\,e^{2\,u}\,u_x^4 + 
          27\,e^u\,u_x^2\,v_x^2 + 
          4\,v_x^4 \right)
          \right.\nn\\
          &&\left.
            - 
     \frac{1 }{864}e^u\,f_{6v}\,
        \left( e^{2\,u}\,u_x^4 + 
          6\,e^u\,u_x^2\,v_x^2 + 
          v_x^4 \right)   -\frac1{216} e^u\,f_{u\,5v}\,
          u_x\,v_x\,
          \left( e^u\,u_x^2 + v_x^2 \right)  
\right\}  +O(\epsilon^6)
\nn
\eeqa
The quasitriviality transformation is generated by the Hamiltonian
\eqa
&&
K=\int dx\, \left\{-\frac{\epsilon}{24}
\left[ 
u\,v_x + 
  (v_x - e^{\frac{u}{2}}\,u_x )\,
  \log(v_x - e^{\frac{u}{2}}\,u_x)  
  + (v_x + e^{\frac{u}{2}}\,u_x )\,
  \log(v_x + e^{\frac{u}{2}}\,u_x) \right]\right.
\nn\\
&&           
+
 \frac{\epsilon^3}{2880 \,\Delta^3}\left[ 2\,e^{3\,u}\,u_x^8\,v_x - 
       13\,e^{3\,u}\,u_x^4\,u_{xx}^2\,
        v_x - 6\,e^{3\,u}\,u_x^2\,
        u_{xx}^3\,v_x - 
       2\,e^{2\,u}\,u_{xx}^3\,v_x^3 + 
       6\,e^{3\,u}\,u_x^5\,u_{xx}\,v_{xx} 
       \right.
       \nn\\
       &&
       \nn\\
       &&\left.
       + 
       6\,e^{3\,u}\,u_x^3\,u_{xx}^2\,
        v_{xx} + 24\,e^{2\,u}\,u_x^3\,
        u_{xx}\,v_x^2\,v_{xx} + 
       18\,e^{2\,u}\,u_x\,u_{xx}^2\,
        v_x^2\,v_{xx} - 
       6\,e^u\,u_x\,u_{xx}\,v_x^4\,
        v_{xx} - 8\,e^{2\,u}\,u_x^4\,v_x\,
        v_{xx}^2 
        \right.
        \nn\\
        &&
        \nn\\
        &&\left.
        - 18\,e^{2\,u}\,u_x^2\,
        u_{xx}\,v_x\,v_{xx}^2 - 
       8\,e^u\,u_x^2\,v_x^3\,
        v_{xx}^2 - 6\,e^u\,u_{xx}\,v_x^3\,
        v_{xx}^2 + 4\,v_x^5\,v_{xx}^2 + 
       2\,e^{2\,u}\,u_x^3\,v_{xx}^3 + 
       6\,e^u\,u_x\,v_x^2\,v_{xx}^3 
\right] 
\nn\\
&&
+
\frac{\epsilon^3}{2880\, \Delta^2}v_x\,
     \left( 2\,e^{2\,u}\,u_x^6 + 
       14\,e^{2\,u}\,u_x^4\,u_{xx} - 
       9\,e^u\,u_x^2\,u_{xx}\,v_x^2 + 
       u_{xx}\,v_x^4 \right) 
\nn\\
&& \left.  
-     
  \frac{\epsilon^3}{1152\, \Delta}e^u\,u_x^4\,v_x
 -\frac{\epsilon^3}{3840}e^{\frac{u}{2}}\,u_x\,
       \left( u_x^2 + 4\,u_{xx} \right) \,
       \log \frac{ v_x+e^{\frac{u}{2}}\,u_x }
         { v_x- e^{\frac{u}{2}}\,u_x   }\right\}
         \nn
 \eeqa
where
$$
\Delta= v_x^2-e^u\,u_x^2 
$$
(cf. \eqref{det-den}).

The main scheme explained above for the case of Boussinesq equation works essentially without changes for the Toda (and for the closely related NLS, see below) hierarchy\footnote{In the recent papers \cite{mm1, mm2} the methods of \cite{du2, dgk} were applied to the description of the critical behaviour of certain particular solutions to the Toda hierarchy and some related systems.}. We will not enter into details here.

\subsection{Nonlinear Schr\"odinger equation}\par

The $D$-operator for focusing NLS is given by
\eqa\label{hf}
&&
f\mapsto h_f=D_{\rm NLS} f = f -\frac{\epsilon^2}{12} \left[ \left( f_{uuu} +\frac3{2u} f_{uu}\right) u_x^2 + 2 f_{uuv} u_x v_x -u f_{uuu} v_x^2\right]
\nn\\
&&
\\
&&
+\epsilon^4 \left\{ \frac1{120} \left[\left( f_{uuuu} +\frac5{2u} f_{uuu}\right) u_{xx}^2 +2 f_{uuuv} u_{xx} v_{xx} -u f_{uuuu} v_{xx}^2\right]\right.
\nn\\
&&
\nn\\
&&
- \frac1{80} f_{uuuu} u_{xx} v_x^2 -\frac1{48 u} f_{uuuv} v_{xx} u_x^2
-\frac1{3456 u^3} \left( 30 f_{uuu} -9 u f_{uuuu} + 12 u^2 f_{5 u} + 4 u^3 f_{6 u}\right) u_x^4
\nn\\
&&
\nn\\
&&
-\frac1{432 u^2} \left( -3 f_{uuuv} + 6 u f_{uuuuv} +2 u^2 f_{5u\, v}\right) u_x^3 v_x +\frac1{288 u} \left( 9 f_{uuuu} + 9 u f_{5 u} + 2 u^2 f_{6u}\right) u_x^2 v_x^2
\nn\\
&&
\nn\\
&&\left.
+\frac1{2160} \left( 9 f_{uuuuv} +10 u f_{5u\, v}\right) u_x v_x^3
-\frac{u}{4320}\left( 18 f_{5u} + 5 u f_{6u}\right) v_x^4\right\}+O(\epsilon^6)
\nn
\eeqa
For the defocusing case one has to replace $u\mapsto -u$. In particular, taking
$$
f=\frac12 (u\, v^2 -u^2)
$$
one recovers the Hamiltonian of the focusing NLS equation
$$
h_f = \frac12 (u\, v^2 -u^2)+\frac{\epsilon^2}{8u} u_x^2.
$$
In this case the infinite series truncates.  Taking $g=-\frac12 v^2 + u (\log u-1)$   one obtains the Hamiltonian of the Toda equation
\eqa
&&
h_g= -\frac12 v^2 + u (\log u-1)-\frac{\epsilon^2}{24 u^2} \left( u_x^2 + 2 u \, v_x^2\right) 
\nn\\
&&
\nn\\
&&
-\epsilon^4 \left( \frac{u_{xx}^2}{240 u^3} +\frac{v_{xx}^2}{60 u^2}
+\frac{u_{xx} v_{x}^2}{40 u^3} - \frac{u_x^4}{144 u^5} -\frac{u_x^2 v_x^2}{24 u^4} +\frac{v_x^4}{360 u^3}\right) +O(\epsilon^6)
\eeqa
written in terms of the function $\phi=\log u$ in the form
$$
\epsilon^2 \phi_{xx} +e^{\phi(s+\epsilon)} -2 e^{\phi(s)} +e^{\phi(s-\epsilon)}=0.
$$
\begin{remark} In \cite{dgk} a somewhat different scaling procedure was used for the description of critical behaviour in the focusing NLS.
\end{remark}

\subsection{Ablowitz - Ladik equation}\par

Ablowitz - Ladik (AL) system \cite{al} is written in the (complexified) form as follows
\eqa\label{abla1}
&&
i\,\dot a_n = -\frac12\,(1-a_n b_n) (a_{n-1}+a_{n+1})+ a_n
\nn\\
&&
\\
&&
i\, \dot b_n = \quad\frac12\,(1-a_n b_n)(b_{n-1} + b_{n+1}) - b_n.
\nn
\eeqa
The reductions
\beq\label{abla2}
b_n =\mp \bar a_n
\eeq
correspond to the focusing/defocusing versions of the AL equation respectively: 
\beq\label{dnls}
i\, \dot a_n +\frac12( a_{n+1} -2 a_n +a_{n-1}) \pm |a_n|^2 \, \frac{a_{n+1}+a_{n-1}}2=0.
\eeq
The latter equation is universally recognized as a discrete version of the NLS equation.

The system \eqref{abla1} is Hamiltonian with respect to the Poisson bracket
\beq\label{abla3}
\{ a_n, b_m\} = i\, (1-a_n b_n) \, \delta_{n,m}, \quad \{ a_n, a_m\} = \{ b_n, b_m\}=0
\eeq
with the Hamiltonian
\beq\label{abla4}
H=\sum_n \frac12(a_n b_{n-1}+b_n a_{n-1}) + \log(1-a_n b_n).
\eeq
Actually the Hamiltonian
\beq\label{abla5}
C=\sum_n \log(1-a_n b_n)
\eeq
generates phase transformations
\eqa
&&
\{ a_n, C\}= -i\, a_n
\nn\\
&&
\{ b_n, C\}=\quad i\, b_n.
\nn
\eeqa
Factoring phase transformations out can be achieved by introducing new variables
\eqa\label{abla6}
&&
w_n =-\log(1-a_n b_n)
\nn\\
&&
\\
&&
v_n = \frac1{2i} \left( \log{\frac{a_n}{a_{n-1}}} -\log{\frac{b_n}{b_{n-1}}} \right).
\nn
\eeqa
Observe that the variables $v_n$, $w_n$ take real values for the focusing/defocusing reductions of AL system. The variable
$w_n$ takes negative values for the focusing case and positive for the defocusing case.
In the new coordinates the dynamics is governed by the Hamiltonian
\beq\label{abla7}
H=\sum_n \sqrt{\left(1-e^{-w_n}\right)\left(1-e^{-w_{n-1}}\right)}\, \cos v_n
\eeq
with the Poisson bracket
\beq\label{abla8}
\{ w_n, v_m\} = \delta_{n,m-1} -\delta_{n,m},\quad \{ w_n,w_m\}=\{ v_n, v_m\}=0.
\eeq
After the interpolation and time rescaling
$$
w_n =w(\epsilon\,n, \epsilon\, t), \quad v_n =v(\epsilon\,n, \epsilon\, t)
$$
one obtains, in the limit $\epsilon\to 0$ the 1st order quasilinear system
\eqa\label{abla9}
&&
w_t =\pal_x\left[\left( e^{-w}-1\right)\, \sin v\right]
\nn\\
&&
\\
&&
v_t\,= \pal_x \left[ e^{-w} \cos v\right].
\nn
\eeqa
We call \eqref{abla9} the {\it long wave limit} of the AL system, in analogy with the long wave limit of Toda lattice.

The characteristic velocities for \eqref{abla9} are
$$
\lambda_\pm =-e^{-w} \left[  \sin v \pm \sqrt{e^{w}-1}\, \cos v\right]
$$
So the long wave limit \eqref{abla9} of the AL system has elliptic type for the focusing case and hyperbolic type for the defocusing case. In the defocusing case the Riemann invariants
of the system \eqref{abla9} are given by the expressions
\beq\label{abla10}
r_\pm = v\pm 2 \arctan \sqrt{e^w-1}.
\eeq

For a later convenience we will replace the dependent variable $w(x)$ by a new variable $u(x)$ by putting
\beq\label{uw1}
u=\frac{\epsilon\, \pal_x}{e^{\epsilon\, \pal_x} -1}\, w
\eeq
(cf. \eqref{uw}). The Poisson brackets of the variables $u$, $v$ take the standard form \eqref{pb}. The Hamiltonian of the interpolated system reads
\eqa\label{abla14}
&&
H_{\rm AL}=\int h_{\rm AL}\, dx = \int  \sqrt{ \left( 1-\exp\{ \frac{1-e^{\epsilon\, \pal_x}}{\epsilon\, \pal_x}\, u\}\right)
\left(1- \exp\{ \frac{e^{-\epsilon\, \pal_x}-1}{\epsilon\, \pal_x}\, u\}\right)}\, \cos v\, dx
\nn\\
&&
\\
&&
h_{\rm AL}=\left( 1-e^{-u}\right)\, \cos v +\frac{\epsilon^2}{24} e^{-u} u_x\,
\frac{(e^u -4) \, \cos v\, u_x +4(e^u-1) \, \sin v \, v_x}{e^u-1}+{\mathcal O}\left(\epsilon^4\right).
\nn
\eeqa
In the long wave limit equations \eqref{abla9} one has just to replace the letter $w$ with the letter $u$.

Let us now establish a connection of the theory of AL equation with Hamiltonian perturbations of the nonlinear wave equation with the potential
\beq\label{abla11}
P(u) =Li_2(e^{-u}).
\eeq 
Recall that the  dilogarithm is obtained by analytic continuation of the function
$$
Li_2(x) =\sum_{k\geq 1} \frac{x^k}{k^2}\quad\mbox{for}\quad |x|<1.
$$
The associated nonlinear wave equation reads
\beq\label{abla12}
u_{tt} =\frac{u_{xx}}{e^u-1}.
\eeq
The linear PDE \eqref{chap1} for the conserved quantities of \eqref{abla12} reads
\beq\label{abla13}
f_{vv}=\left( e^u-1\right)\, f_{uu}.
\eeq

\begin{lemma} The function $f(u,v)$ is a conserved quantity for the long wave limit \eqref{abla9} of AL system {\rm iff} it satisfies \eqref{abla13}. 
\end{lemma}

\pf is given by a straightforward calculation. \epf

We can now apply the above arguments to constructing the $D$-operator for the AL hierarchy.

\begin{theorem} For any ${\mathcal C}^\infty$ solution $f=f(u,v)$ to \eqref{abla13} there exists a unique Hamiltonian 
\eqa\label{abla15}
&&
H_f = H_f^{[0]} +\sum_{k\geq 1} \epsilon^k H_f^{[k]} = \int D_{\rm AL}f\, dx
\nn\\
&&
\\
&&
H_f^{[k]}=\int h_f^{[k]}\, dx, \quad h_f^{[0]}=f
\nn
\eeqa
commuting with the {\rm AL} Hamiltonian \eqref{abla14}:
$$
\int \left[ \frac{\delta H}{\delta u(x)} \pal_x  \frac{\delta H_{\rm AL}}{\delta v(x)}+\frac{\delta H}{\delta v(x)} \pal_x  \frac{\delta H_{\rm AL}}{\delta u(x)}\right]\, dx=0.
$$
\end{theorem}

Explicitly
\eqa
&&
D_{\rm AL}f= f+
\epsilon^2\left[
 -\frac1{24} \,\frac{2(1-e^u)\, f_{uvv} +(e^u-2)\, f_{vv}}{(e^u-1)^2}\, u_x^2 \right.
\nn\\
&&\left.
\quad\quad\quad+\frac16 \, \frac{f_{vvv}}{e^u-1}\, u_x v_x +\frac1{12}\, \frac{(e^u-1)\, f_{uvv}-f_{vv}}{e^u-1}\, v_x^2\right]+{\mathcal O}\left( \epsilon^4\right).
\nn
\eeqa

The basis of linearly independent solutions to the linear PDE \eqref{abla13} can be obtained in the form
\beq\label{ges}
\begin{array}{l}f_n ^+=\cos n v \, e^{-n\, u} F(-n, -n; 1-2n; e^u) \\
 \\
  f_n ^-=\sin n v \, e^{-n\, u} F(-n, -n; 1-2n; e^u)\end{array},\quad n\in \mathbb Z, \quad n\neq 0.
\eeq
Here $F(\alpha,\beta;\gamma;x)$ is the Gauss hypergeometric function. For positive $n$ the hypergeometric factor is a polynomial in $e^u$ of the degree $n$; these conserved quantities correspond to the first integrals of AL hierarchy constructed in \cite{al}. The second part of the solutions to \eqref{abla13}, the one for $n<0$, corresponds to the new integrals of the AL hierarchy found recently in \cite{ges}.

For example, choosing 
$$
f=f_1^+=(1-e^{-u})\, \cos v
$$
one obtains
$$
h_f = h_{\rm AL}.
$$
The choice
$$
f=f_1^-=(1-e^{-u})\, \sin v
$$
yields the $\epsilon$ expansion 
$$
h_f = (1-e^{-u})\, \sin v +\frac{\epsilon^2}{24} \, e^{-u} u_x\, \frac{(e^u-4)\, \sin v\, u_x -4 (e^u -1)\, \cos v\, v_x}{e^u-1} +{\mathcal O}\left( \epsilon^4\right)
$$
of the Hamiltonian
$$
H=\frac1{2i} \sum a_n b_{n-1} - b_n a_{n-1}
$$
of the discrete modified KdV equation
\eqa
&&
\dot a_n =\frac12\,(1-a_n b_n) (a_{n+1}-a_{n-1})
\nn\\
&&
\dot b_n =\frac12\,(1-a_n b_n) (b_{n+1}-b_{n-1}).
\nn
\eeqa
Observe that the Hamiltonian of the perturbed wave equation \eqref{abla12} reads
\beq\label{abla16}
H_{\rm pert}=\int \left[ \frac12\, v^2 +Li_2(e^{-u}) + \frac{\epsilon^2}{24}
\frac{u_x^2 -(e^u-1)\, (u_x^2 + 2 v_x^2)}{(e^u-1)^2}\right]\, dx+{\mathcal O}\left( \epsilon^4\right).
\eeq
This perturbation is not equivalent to the one given in \eqref{hami4} coming from the generalized FPU system with the potential \eqref{abla11}.

\begin{remark} The recurrent procedure of constructing the conserved quantities of the long wave limit of the Ablowitz -- Ladik hierarchy can be described in terms of the following solution to WDVV
\beq\label{polylog3}
F(u,v)=\frac12 u\, v^2 -Li_3 \left( e^{-u}\right)
\eeq
(see Remark \ref{rem_fro} above).
\end{remark}

%\newpage


\begin{thebibliography}{99}

\bibitem{al} M.J.Ablowitz, J.F.Ladik, Nonlinear differential-difference equations and Fourier analysis, {\it J. Math. Phys.} {\bf 17} (1976) 1011-1018. 

\bibitem{ar} V.I.Arnold,V.V.Goryunov, O.V.Lyashko,V.A.Vasil'ev,  {\it Singularity Theory. I.}  Dynamical systems. VI, Encyclopaedia Math. Sci. {\bf 6}, Springer, Berlin, 1993.

\bibitem{bgi} V.A.Baikov, R.K. Gazizov, N.Kh. Ibragimov, Approximate symmetries
and formal linearization, {\it PMTF} {\bf 2} (1989) 40-49. (In Russian)

\bibitem{bo} P.Boutroux, Recherches sur les transcendants de M. Painlev\'e et l'\'etude asymptotique des \'equations diff\'erentielles du second ordre. {\it Ann. \'Ecole Norm} {\bf 30} (1913) 265 - 375.

\bibitem{bressan} A.Bressan, One dimensional hyperbolic systems of conservation laws. Current developments in mathematics, 2002, 1--37, Int. Press, Somerville, MA, 2003.

\bibitem{bmp} \'E.Br\'ezin,  E.Marinari, G.Parisi, 
A nonperturbative ambiguity free solution of a string model. 
{\it Phys. Lett.} {\bf B 242} (1990) 35--38.

\bibitem{cdz} G.Carlet, B.Dubrovin, Y.Zhang, The Extended Toda Hierarchy,  {\it Moscow Math. J.} {\bf 4} (2004) 313-332.

\bibitem{cg} T.Claeys, T.Grava, Universality of the break-up profile for the KdV equation in the small dispersion limit using the Riemann-Hilbert approach, arXiv:0801.2326.

\bibitem{cl1} T.Claeys, M.Vanlessen, The existence of a real pole-free solution of the fourth order analogue of the Painlev\'e I equation. {\it Nonlinearity} {\bf 20} (2007), no. 5, 1163--1184. 

\bibitem{cl2} T. Claeys, M. Vanlessen, Universality of a double scaling limit near singular edge points in random matrix models, {\it Comm. Math. Phys.} {\bf 273} (2007), no. 2, 499--532.

\bibitem{dt} P.Dedecker and W.M.Tulczyjev, {\it Spectral Sequences and the Inverse
Problem of the Calculus of Variations}, Lecture Notes in Math. {\bf 836} (1980)
498-503.

\bibitem{magri} L. Degiovanni, F. Magri, V. Sciacca, On deformation of Poisson
manifolds of hydrodynamic type, {\it Comm. Math. Phys.} {\bf 253} (2005), no. 1, 1--24.

\bibitem{dz} P.Deift, X.Zhou, A steepest descent method for oscillatory Riemann-Hilbert problems. Asymptotics for the MKdV equation, {\it Ann. of Math.} (2) {\bf 137} (1993), 295--368.


\bibitem{npb} B. Dubrovin, Integrable systems in topological
field theory, {\it Nucl. Phys.} {\bf B379} (1992), 627--689.

\bibitem{D1} B. Dubrovin, Geometry of 2D topological field theories, in {\sl Integrable systems and quantum
groups} (Montecatini Terme, 1993), 120 - 348. Lecture Notes in Mathematics, {\bf 1620}. Springer, Berlin, 1996.

\bibitem{icm} B. Dubrovin, Geometry and analytic theory of Frobenius manifolds,
Proceedings of ICM98, Vol. 2, 315-326.

\bibitem{du2} B.Dubrovin, On Hamiltonian perturbations of hyperbolic systems of conservation laws, II: universality of critical behaviour, 
{\it Comm. Math. Phys.} {\bf 267} (2006) 117 - 139.


\bibitem{dgk} B.Dubrovin, T.Grava, C.Klein, On universality of critical behaviour in the focusing nonlinear Schr\"odinger equation, elliptic umbilic catastrophe and the {\it tritronquŽe} solution to the Painlev\'e-I equation, arXiv:0704.0501, to appear in {\it J. Nonlinear Sci.}

\bibitem{DLZ} B. Dubrovin, S.-Q. Liu, Y. Zhang, On Hamiltonian perturbations of hyperbolic systems
of conservation laws, I: Quasi-triviality of bi-Hamiltonian perturbations. {\it Comm. Pure Appl. Math.} {\bf 59} (2006),
no. 4, 559 - 615.


\bibitem{dn83} B. Dubrovin, S.P. Novikov, The Hamiltonian formalism of one-dimensional systems of hydrodynamic
type and the Bogolyubov -- Whitham averaging method,
{\it Soviet Math. Dokl.} {\bf 270:4} (1983), 665-669.
\bibitem{dn84} B. Dubrovin, S.P. Novikov,
On Poisson brackets of hydrodynamic type,
{\it Soviet Math. Dokl.} {\bf 279:2} (1984), 294-297.
\bibitem{dn89} B. Dubrovin, S.P. Novikov, Hydrodynamics of weakly deformed
soliton lattices. Differential geometry and Hamiltonian theory. {\it Uspekhi Mat.
Nauk} {\bf 44} (1989) 29-98. English translation in {\it Russ. Math. Surveys} {\bf 44}
(1989) 35-124.


\bibitem{DZ} B.Dubrovin, Y.Zhang, Normal forms of hierarchies of integrable PDEs, Frobenius manifolds
and Gromov - Witten
invariants, ArXiv:math.DG/0108160, 2001.


\bibitem{ges} F.Gesztesy, H.Holden, J.Michor, G.Teschl, Local conservation laws and the Hamiltonian formalism for the Ablowitz-Ladik hierarchy, ArXiv:0711.1644

\bibitem{getzler} E. Getzler, A Darboux theorem for Hamiltonian operators
in the formal calculus of variations, {\it Duke Math. J.} {\bf 111} (2002), 535--560.

\bibitem{gk1} T.Grava, C.Klein, Numerical solution of the small dispersion limit of Korteweg de Vries and Whitham equations. {\it Comm. Pure Appl. Math.} {\bf 60} (2007), no. 11, 1623--1664.

\bibitem{gk} T.Grava, C.Klein, Numerical study of a multiscale expansion of KdV and Camassa-Holm equation. ArXiv:math-ph/0702038.

\bibitem{GN} P.Grinevich, S.P. Novikov,  String equation. II. Physical solution. (Russian)  
{\it Algebra i Analiz}  {\bf 6}  (1994),  no. 3, 118--140;  translation in  {\it St. Petersburg Math. J.}  {\bf 6}  (1995) 553--574.

\bibitem{gp} A.Gurevich, L.Pitaevski, Nonstationary structure of a collisionless
shock wave, {\it Sov. Phys. JETP Lett.} {\bf 38} (1974) 291--297.

\bibitem{hk1} A.Henrici, T.Kappeler, Normal forms for odd periodic FPU chains, arXiv:nlin/0611063.

\bibitem{hk2} A.Henrici, T.Kappeler, Resonant normal form for even periodic FPU chains, arXiv:0709.2624.

\bibitem{hl} T.Y.Hou,  P.D.Lax, 
Dispersive approximations in fluid dynamics. 
{\it Comm. Pure Appl. Math.} {\bf 44} (1991) 1--40.

\bibitem{kapaev} A.A.Kapaev,  
Weakly nonlinear solutions of the equation ${\rm P}\sp 2\sb 1$,
{\it Zap. Nauchn. Sem. Leningrad. Otdel. Mat. Inst. Steklov. (LOMI)} {\bf 187} (1991), Differentsialnaya Geom. Gruppy Li i Mekh. {\bf 12}, 88--109, 172--173, 175; translation in {\it J. Math. Sci.} {\bf 73} (1995), no. 4, 468--481.

\bibitem{ka} A.Kapaev, Quasi-linear Stokes phenomenon for the Painlev\'e first equation. {\it J. Phys.  A: Math. Gen.} {\bf 37} (2004) 11149--11167.

\bibitem{km} Y. Kodama, A. Mikhailov, 
Obstacles to asymptotic integrability, 
Algebraic aspects of integrable systems, 173--204,
{\it Progr. Nonlinear Differential Equations Appl.}, {\bf 26},
Birkh\"auser, Boston, MA, 1997.

\bibitem{ks} V.Kudashev, B.Suleimanov, A soft mechanism for the generation of dissipationless shock waves, {\it Phys. Lett.} {\bf A 221} (1996) 204--208.

\bibitem{lax} P.Lax, Periodic solutions of the KdV equations. {\it Nonlinear Wave Motion}, pp. 85--96. {\it Lectures in Appl. Math.}, {\bf 15}, Amer. Math. Soc., Providence, R.I., 1974.

\bibitem{ll} P.Lax, D.Levermore,
The small dispersion limit of the Korteweg-de Vries equation. I, II, III.
{\it Comm. Pure Appl. Math.} {\bf 36} (1983) 253--290, 571--593, 809--829.
\bibitem{llv} P. D.Lax, C. D.Levermore, S.Venakides,
The generation and propagation of oscillations in dispersive initial value problems and their limiting behavior. In: 
{\it Important Developments in Soliton Theory}, 205--241, 
Springer Ser. Nonlinear Dynam., 
Springer, Berlin, 1993. 

\bibitem{LWZ} S.-Q. Liu, C.-Z.Wu, Y. Zhang, On properties of Hamiltonian structures for a class of evolutionary PDEs, ArXiv:0711.2599.

\bibitem{LZ2} S.-Q. Liu, Y. Zhang, On Quasitriviality and Integrability of a Class of Scalar Evolutionary PDEs,
{\it J. Geom. Phys.} {\bf 57} (2006) 101-119.

\bibitem{mm1} L.Mart\'\i nez Alonso, E.Medina, Regularization of Hele-Shaw flows, multiscaling expansions and the Painlev\'e I equation, arXiv:0710.3731.

\bibitem{mm2} L.Mart\'\i nez Alonso, E.Medina,
The double scaling limit method in the Toda hierarchy, arXiv:0804.3498.

\bibitem{n74} S.P.Novikov, A periodic problem for the Korteweg-de Vries equation. I. {\it Funkcional. Anal. i Prilo\v zen.} {\bf 8 }(1974), no. 3, 54--66.

\bibitem{serre} D.Serre, {\it Syst\`emes de lois de conservation I : hyperbolicit\'e, entropies, ondes de choc}; {\it Syst\`emes de lois de conservation II: structures g\'eom\'etriques, oscillation et probl\`emes mixtes},
Paris Diderot Editeur 1996.

\bibitem{str} I.A.B.Strachan, Deformations of the Monge/Riemann hierarchy and approximately integrable systems,
 {\it J. Math. Phys.}  {\bf 44}  (2003) 251--262.

\bibitem{th} R.Thom, {\it Structural Stability and Morphogenesis: An Outline of a General Theory of Models.} Reading, MA: Addison-Wesley, 1989.

\bibitem{tsarev} S. Tsarev, The geometry of Hamiltonian systems
of hydrodynamic type. The generalized hodograph method,
{\it Math. USSR Izv.} {\bf 37} (1991), 397--419.

\bibitem{w} H.Whitney, On singularities of mappings of euclidean spaces. I. Mappings of the plane into the plane. 
{\it Ann. of Math.} (2) {\bf 62} (1955), 374--410. 

\bibitem{zk} N.Zabusky, M.Kruskal, {\it Phys. Rev. Lett.} {\bf 15} (1965) 2403.


\end{thebibliography}
\end{document}